\documentclass[draft,onecolumn]{IEEEtran}
\ifCLASSINFOpdf
   \usepackage[pdftex]{graphicx}
    \usepackage{multicol}
    \usepackage{multirow}
    \usepackage{stfloats}
\else
\fi

\newcommand{\nn}{\nonumber}
\newcommand{\bmat}{\left[\begin{array}}
\newcommand{\emat}{\end{array}\right]}
\usepackage[table]{xcolor}
\usepackage{graphicx}
\usepackage{epsfig}
 \usepackage{amsmath}
\usepackage{amssymb}
\usepackage{amsfonts}
\usepackage{amsthm}

 \usepackage{epstopdf}
\usepackage{algorithm}
\usepackage{algpseudocode}
\usepackage{cite}
\usepackage{subfigure}
\usepackage{bbm}
\definecolor{mygreen}{RGB}{34,139,34}

\newcommand{\tr}{\operatorname{tr}}
\newcommand{\rank}{\operatorname{rank}}
\newcommand{\Rs}{\mathbb R}

\newcommand{\Es}{\mathbb{E}}
\newtheorem{remark}{Remark}
\newtheorem{lemma}{Lemma}
\newtheorem{corollary}{Corollary}
\newtheorem{proposition}{Proposition}
\newtheorem{theorem}{Theorem}
\newcommand{\al}[1]{\begin{align} #1 \end{align}}
\begin{document}

\title{Robust Kalman Filtering Under Model Uncertainty:\\ the Case of Degenerate Densities
}

\author{Shenglun Yi and Mattia Zorzi, {\em Senior Member IEEE} \thanks{S. Yi  is with the School of Automation, Beijing Institute of Technology, Beijing 100081, China; M. Zorzi is with the Department of Information Engineering, University of Padova, Via Gradenigo 6/B, 35131 Padova, Italy. Emails: {\tt\small 3120185460@bit.edu.cn}, {\tt\small zorzimat@dei.unipd.it}}
\thanks{{ This work is partially supported by the SID project ``A Multidimensional and Multivariate Moment Problem Theory for Target Parameter Estimation in Automotive Radars'' (ZORZ\_SID19\_01) funded by the Department of Information Engineering of the University of Padova, Italy.}
}
}

\maketitle

\begin{abstract}  We consider a robust state space filtering problem in the case that the transition probability density is unknown and possibly degenerate. The resulting robust filter has a Kalman-like structure and solves a minimax game: the nature selects the least favorable model in a prescribed ambiguity set which also contains non-Gaussian probability densities, while the other player designs the optimum filter for the least favorable model. It turns out that the resulting robust filter is characterized by a Riccati-like iteration evolving on the cone of the positive semidefinite matrices. Moreover, we study the convergence of such iteration in the case that the nominal model is with constant parameters on the basis of the contraction analysis in the same spirit of Bougerol. Finally, some numerical examples show that the proposed filter outperforms the standard Kalman filter.
\end{abstract}

\begin{IEEEkeywords}
 Robust Kalman filtering, minimax problem, low-rank filtering, least favorable model, contraction analysis.
\end{IEEEkeywords}

\IEEEpeerreviewmaketitle
\section{Introduction}
As a classic estimation algorithm, the standard Kalman filter is extensively used. However, it may perform poorly in the case that the actual model does not coincide with the nominal one. For such a reason, manifold robust versions of the Kalman filter have been considered, see for instance \cite{HASSIBI_SAYED_KAILATH_BOOK,speyer_2008,GHAOUI_CALAFIORE_2001,kim2020robust,Li2020}.

In particular, risk sensitive Kalman filters \cite{RISK_WHITTLE_1980,RISK_PROP_BANAVAR_SPEIER_1998,LEVY_ZORZI_RISK_CONTRACTION,huang2018distributed,OPTIMAL_SPEYER_FAN_BANAVAR_1992} have been proposed in order to address model uncertainty. The latter consider an exponential quadratic loss function rather than the standard quadratic loss function. This means that large errors are severely penalized according to the so called risk sensitivity parameter. Hereafter, these robust Kalman filters have been proved to be equivalent to solve a minimax problem, \cite{boel2002robustness,HANSEN_SARGENT_2005,YOON_2004,HANSEN_SARGENT_2007}. More precisely, there are two players. One player, say nature, selects the least favorable model in a prescribed ambiguity set which is a ball about the nominal model and whose radius reflects the amount of uncertainty in respect to the nominal model. The other player designs the optimum filter according to the least favorable model.

Recently, instead of concentrating the entire model uncertainty in a single relative entropy constraint, a new paradigm of risk sensitive filters has been proposed in \cite{ROBUST_STATE_SPACE_LEVY_NIKOUKHAH_2013,zorzi2018robust,abadeh2018wasserstein,zorzi2019distributed,RKDISTR_opt,zenere2018coupling}. The latter characterizes the uncertainty using separate constraints to each time increment of the model. In other words, the ambiguity set is specified at each time step by forming a ball, in the Kullback-Leibler (KL) topology, about the nominal model, \cite{ROBUST_STATE_SPACE_LEVY_NIKOUKHAH_2013,robustleastsquaresestimation}. It is worth noting this ball can be defined by using also the Tau-divergence family,  \cite{STATETAU_2017,OPTIMALITY_ZORZI}. These filters, however, are well defined only in the case that the nominal and the actual transition probability density functions are non-degenerate. This guarantees that the corresponding distorted Riccati iteration evolves on the cone of the positive definite matrices.

Unfortunately, in many practical applications, like weather forecasting and oceanography, the standard Kalman filter fails to work. More precisely, the Riccati iteration produces numerical covariance matrices which are indefinite because of their large dimension. This issue is typically addressed by resorting to a low-rank Kalman algorithm \cite{dee1991simplification}.

The contribution of this paper is to extend the robust paradigm in \cite{ROBUST_STATE_SPACE_LEVY_NIKOUKHAH_2013} to the case in which the transition probability density is possibly degenerate. Some preliminary results can be found in \cite{yi2020low}. Within our framework, degenerate Gaussian probability densities could be also involved in the dynamic game. Accordingly, the resulting robust Kalman filter corresponds to a low-rank risk sensitive Riccati iteration. Although low-rank and distorted Riccati iterations have been widely studied in the literature, e.g. \cite{bonnabel2013geometry,Ferrante20141176,ferrante2013generalised,Baggio-Ferrante-TAC-19,Baggio-Ferrante-TAC-16}, our iteration appears to be new.  Then, we also derive the least favorable model over a finite simulation horizon. Last but not least, we study the convergence of the distorted Riccati iteration in the case that the nominal model has constant parameters by means of the contraction analysis, \cite{BOUGEROL_1993}. It turns out that  the convergence is guaranteed if the nominal model is stabilizable, the reachable subspace is observable and the ambiguity set is ``small''.

The outline of the paper is as follows. Section \ref{sec_2} discusses the low-rank robust static estimation problem where the ambiguity set is characterized by a relative entropy constraint and possibly contains degenerate densities. The robust Kalman filtering problem is presented in Section \ref{sec_3}. The latter is then reformulated as a static minimax game in Section \ref{sec_4}. In Section \ref{sec_5}, we derive the corresponding least favorable model. Section \ref{sec_6} regards the convergence of the proposed low-rank robust Kalman filter in the case of constant parameters. In Section \ref{sec_7} some numerical examples are provided. Finally, we draw the conclusions in Section \ref{sec_8}.

{\em Notation:} The image, the kernel and the trace of matrix $K$  are denoted by $\mathrm{Im}(K)$, $\mathrm{ker}(K)$ and $\mathrm{tr}(K)$, respectively.  Given a symmetric matrix $K$:  $K>0$ ($K\geq 0$) means that $K$ is positive (semi)definite; $\sigma_{max}(K)$ is the maximum eigenvalue of $K$; $K^+$ and $\det^+(K)$ denote the pseudo inverse and the pseudo determinant of $K$, respectively. The Kronecker product between two matrices $K$ and $V$ is denoted by  $K\otimes V$. $x\sim \mathcal N(m,K)$ means that $x$ is a Gaussian random variable with mean $m$ and covariance matrix $K$. $\mathcal Q^n$ is the vector space of symmetric matrices of dimension $n$; $\mathcal Q_+^n$ denotes the cone of positive definite symmetric matrices of dimension $n$, while $ \overline{\mathcal Q_+^n}$ denotes its closure. The diagonal matrix whose elements in the main diagonal are $a_1,a_2,\ldots , a_n$ is denoted by $\mathrm{diag}(a_1,a_2, \ldots a_n)$; $\mathrm{Tp}(A_1, A_2,\ldots ,A_n)$ denotes the block upper triangular Toeplitz matrix whose first block row is $[\,A_1 \; A_2\; \ldots \; A_n \, ]$.
 \section{Low-rank robust static estimation}\label{sec_2}
Consider the problem to estimate a random vector $x\in\mathbb R^n$ from an observation vector $y\in\mathbb R^p$. We assume that the nominal probability density function of $z:=[\, x^{\top}\; y^{\top} \,]^{\top}$ is $f(z) \sim \mathcal{N}\left(m_{z}, K_{z}\right)$ where the mean vector $m_z\in\mathbb R^{n+p}$ and the covariance matrix $K_z\in \overline{\mathcal Q^{n+p}_+}$ are such that
$$
m_{z}=\left[\begin{array}{c}{m_{x}} \\ {m_{y}}\end{array}\right], \quad K_{z}=\left[\begin{array}{cc}{K_{x}} & {K_{x y}} \\ {K_{y x}} & {K_{y}}\end{array}\right].
$$
Moreover, we assume that $K_z$ is possibly a singular matrix and such that $\mathrm{rank}(K_z)=r+p$ with $r\leq n$ and $K_y>0$. Therefore, $f(z)$ is possibly a degenerate density whose support is the $r+p$-dimensional affine subspace
$$
\mathcal{A}=\left\{m_{z}+v, \quad v \in \mathrm {Im} \left(K_{z}\right)\right\}
$$
and
\begin{equation} \label{pdf_f}
\begin{aligned} f(z) =&\left[(2 \pi)^{r+p} \operatorname{det}^{+}\left(K_{z}\right)\right]^{-1 / 2}  \\
& \times  \exp \left[-\frac{1}{2}\left(z-m_{z}\right)^{\top} K_{z}^{+}\left(z-m_{z}\right)\right]. \end{aligned}
\end{equation}
Let $\tilde f(z) \sim \mathcal{N}(\tilde m_{z}, \tilde K_{z})$ denote the actual probability density function of $z$ and we assume that $\mathrm{rank}(\tilde K_z)=r+p$. Accordingly,
\begin{equation}\label{pdf_tildef}
\begin{aligned} \tilde f(z) =&\left[(2 \pi)^{{r+p}} \operatorname{det}^{+}\left(\tilde K_{z}\right)\right]^{-1 / 2}   \\
& \times \exp \left[-\frac{1}{2}\left(z-\tilde m_{z}\right)^{ \top} \tilde K_{z}^{+}\left(z-\tilde m_{z}\right)\right] \end{aligned}
\end{equation}
with support
$\mathcal{\tilde A}= \{\tilde m_{z}+v, \quad v \in \mathrm {Im}(\tilde{ K}_{z}) \}.$ We use the KL-divergence to measure the deviation between the nominal probability density function $f(z)$ and the actual one $\tilde f(z)$. Since the KL-divergence is not able to measure the ``deterministic'' deviations, we have to assume that the two probability density functions have the same support, i.e. $\mathcal{A}=\mathcal{\tilde A}$. In other words, we have to impose that:
\begin{equation} \mathrm {Im}\left(K_{z}\right)=\mathrm {Im} (\tilde{ K}_{z}), \; \; \Delta m_z \in \mathrm {Im}\left(K_{z}\right) \label{KL}\end{equation}
where $\Delta m_z=\tilde m_z-m_z$. Hence, under  assumption (\ref{KL}), the KL-divergence between $\tilde f(z)$ and $f(z)$ is defined as
\begin{equation} \label{def_DL_ddeg}
D (\tilde{f}, f )=\int_\mathcal{A} \ln \left(\frac{\tilde{f}(z)}{f(z)}\right) \tilde{f}(z) d z.
\end{equation}
Then, if we substitute (\ref{pdf_f}) and (\ref{pdf_tildef}) in (\ref{def_DL_ddeg}), we obtain
\begin{equation*}
\begin{aligned} D(\tilde f, {f})=& \frac{1}{2}\left[\Delta{m}_{z}^{\top}K_{z}^{+} \Delta m_z+\ln \operatorname{det}^{+} (K_{z})\right.\\ &\left.-\ln \operatorname{det}^{+} (\tilde{K}_{z})+\tr\left(K_{z}^{+} \tilde{K}_{z}\right)-(r+p)\right] .\end{aligned}
\end{equation*}

\begin{lemma}
\label{lemma1}
Let $f(z) \sim \mathcal{N}\left(m_{z}, K_{z}\right)$ and $\tilde f(z) \sim \mathcal{N}(\tilde m_{z}, \tilde K_{z})$ be Gaussian and possibly degenerate probability density functions with the same $r+p$-dimensional support $\mathcal{A}$. Let \begin{align*}\mathcal U&=\{\tilde m_z \in\mathbb R^{n+p} \hbox{ s.t. } \tilde m_z-m_z\in \mathrm{Im}(K_z)\}\\
\mathcal{V}&={\{\tilde K_{z} \in \mathcal{Q}^{n+p} ~\text {s.t.}~  \operatorname{Im}(K_{z} )=\operatorname{Im} (\tilde{K}_{z} )}\}. \end{align*} Then, $D(\tilde f,f)$ is strictly convex with respect to $\tilde m_z\in\mathcal U$ and $\tilde K_z\in\mathcal V \cap \overline{\mathcal Q_+^{n+p}}$. Moreover, $D(\tilde f, f) \geq 0$ and equality  holds if and only if $f=\tilde f$.
\end{lemma}


\IEEEproof
Let $U_{\mathfrak{r}}^\top$ be a matrix whose columns form an orthonormal basis for $\mathrm {Im}(K_z)$. Moreover, we define $  K_z^{\mathfrak{r}}= U_{\mathfrak{r}}  K_zU_{\mathfrak{r}}^\top$ and $\tilde K_z^{\mathfrak{r}}= U_{\mathfrak{r}} \tilde K_zU_{\mathfrak{r}}^\top$.  Since $f(z)$ and $\tilde f(z)$ have the same support, then $\tilde K_z$ and $\tilde m_z$ are such that $\mathrm {Im}(\tilde K_z)=\mathrm {Im}( K_z)$ and $\Delta m_z \in \mathrm {Im}(K_z)$.  Thus, $D(\tilde f,f)$ is strictly convex in $\tilde K_z \in\mathcal V\cap  \overline{\mathcal Q_+^{n+p}}$ if and only if it is  strictly convex in $\tilde K_z^{\mathfrak{r}} \in \mathcal Q_+^{r+p}$.  Then, it is not difficult to see that
\begin{equation*}
\begin{aligned} D(\tilde f, {f})=& \frac{1}{2}\left[\Delta m^{\top}_z U^{\top}_{\mathfrak{r}} (K^{\mathfrak{r}}_{z})^{-1} U_{\mathfrak{r}} \Delta m_z +\tr((K^{\mathfrak{r}}_{z})^{-1} \tilde{K}^{\mathfrak{r}}_{z}) \right.\\ &\left.-\ln \operatorname{det} ((K^{\mathfrak{r}}_{z})^{-1} \tilde{K}^{\mathfrak{r}}_{z})-(r+p)\right]\end{aligned}
\end{equation*}
which is strictly convex in $K^{\mathfrak{r}}_{z} \in  \mathcal Q_+^{n+p}$, see \cite{COVER_THOMAS}. Hence, we proved the strict convexity of $D(\tilde f,f)$ in $\tilde K_z \in \mathcal V\cap \overline  {\mathcal Q_+^{n+p}}$. Using  similar reasonings, we can conclude that $D(\tilde f, {f})$ is strictly convex in $\tilde m_z\in \mathcal U$.
Finally, the unique minimum of $D(\tilde f,f)$ with respect to $\tilde f$ is given by the stationary conditions $\tilde m_z=  m_z$ and $\tilde K _{z}={K}_{z}$, i.e. $\tilde f=f$. Since $D(f,f)=0$, we conclude that  $D(\tilde f,f)\geq 0$ and equality holds if and only if $\tilde f=f$.
\qed \\

In what follows, we assume that $f$ is known while $\tilde f$ is not and both have the same support. Then, we assume that $\tilde f$ belongs to the ambiguity set, i.e. a ``ball'':
$$
\mathcal{B}=\{\tilde{f} \sim \mathcal N(\tilde m_z,\tilde K_z)\hbox{ s.t. } D(\tilde{f}, f) \leq c\}
$$
where $c>0$, hereafter called tolerance, represents the radius of this ball.
It is worth noting that $f$ is usually estimated from  measured data. More precisely,  we fix a parametric and Gaussian model class $\mathcal M$, then we select $f\in\mathcal M$ according to the maximum likelihood principle.  Thus, when  the length of the data is sufficiently large, we have
\al{f \approx \underset{f\in\mathcal M}{\mathrm{argmin}}\, D(\tilde f,f)\nn}
under standard hypotheses. Therefore, the uncertainty on $f$ is naturally defined by $\mathcal B$, i.e. the actual model $\tilde f$ satisfies the constraint $\tilde f\in\mathcal B$ with $c=D(\tilde f,f)$. Finally, an estimate of $c$ is given by $\hat c=D(\check f,f )$ where $\check f$ is estimated from measured data using a model class $\check{\mathcal  M}$ sufficiently ``large'', i.e. it contains many candidate models having diverse features.

In view of Lemma \ref{lemma1}, $\mathcal B$ is a convex set. We consider the robust estimator $\hat x=g^0(y)$ solving the following minimax problem
\begin{equation} \label{robust_p}
(\tilde f^0, g^0)=\operatorname{arg}\min _{g \in \mathcal{G}} \max _{\tilde{f} \in \mathcal{B}} J(\tilde{f}, g)
\end{equation}
where$$
\begin{aligned} J(\tilde{f}, g) &=\frac{1}{2} E_{\tilde{f}}\left[\| x-g(y )\|^{2}\right] \\ &=\frac{1}{2} \int_\mathcal{A}\| x-g(y) \|^{2} \tilde{f}(z) d z. \end{aligned}
$$$\mathcal{G}$ is the set of estimators for which  $E_{\tilde{f}}\left[\|x-g(y)\|^{2}\right]$ is bounded with respect to all the Gaussian densities in $\mathcal B$.

\begin{theorem} \label{theo1} Let $f$ be a Gaussian (possibly degenerate) density defined as in (\ref{pdf_f}) with $K_y>0$. Then, the least favorable Gaussian density $\tilde f^0$ is with mean vector and covariance matrix
\al{\label{def_m_K_tilde}
\tilde{m}_{z}^0=m_{z}=\left[\begin{array}{c}{m_{x}} \\ {m_{y}}\end{array}\right], \quad \tilde{K}_{z}^0=\left[\begin{array}{cc}{\tilde{K}_{x}} & {K_{x y}} \\ {K_{y x}} & {K_{y}}\end{array}\right]}
so that, only the covariance of $x$ is perturbed.
Then, the Bayesian estimator
\begin{equation*}
g^{0}(y)=G_{0}\left(y-m_{y}\right)+m_{x},
\end{equation*}
with $G_0=K_{x y} K_{y}^{-1}$, solves the robust estimation problem.
The nominal posterior covariance matrix of $x$ given $y$ is given by \begin{align} \label{def_P_stat}P&:=K_{x}-K_{x y} K_{y}^{-1} K_{y x}.\end{align} while the least favorable one is:
\begin{equation*} \tilde P=(P^+-\lambda ^{-1}H^{ \top} H)^+\end{equation*}
where $H^\top\in\Rs^{n\times r}$ is a matrix whose columns form an orthonormal basis for $ \mathrm {Im} (P)$.
Moreover, there exists a unique Lagrange multiplier $\lambda >\sigma_{max}(P)>0$ such that $c=D(\tilde f^0, f)$.

\end{theorem}

\IEEEproof
The saddle point of $J$ must satisfy the conditions
\begin{equation} \label{ineq_saddle}
J (\tilde{f}, g^{0} ) \leq J (\tilde{f}^{0}, g^{0} ) \leq J(\tilde{f}^{0}, g )
\end{equation}
for all $\tilde f \in \mathcal{B}$ and $g \in \mathcal{G}$.
The second inequality in (\ref{ineq_saddle}) is based on the fact that the Bayesian estimator $g^0$ minimizes $J(\tilde f^0,\cdot)$.  Therefore, it remains to prove the first inequality in (\ref{ineq_saddle}).

Notice that the minimizer of $J(\tilde f, \cdot)$ is $$g^{\star}(y)=\tilde K_{xy}\tilde K_y^{-1}(y-\tilde m_y)+\tilde m_x$$ where
$$
\tilde m_{z}:=\left[\begin{array}{c}{\tilde m_{x}} \\ {\tilde m_{y}}\end{array}\right], \quad \tilde K_{z}:=\left[\begin{array}{cc}{\tilde K_{x}} & {\tilde K_{x y}} \\ {\tilde K_{y x}} & {\tilde K_{y}}\end{array}\right].
$$ Moreover, $J(\tilde f,g^\star)=1/2\tr(\tilde P)$ where $\tilde P:=\tilde K_{x}-\tilde K_{x y} \tilde K_{y}^{-1} \tilde K_{y x}$. Since $\mathrm{Im}(\tilde K_z)=\mathrm{Im}(K_z)$, then $\mathrm{Im}(\tilde P)=\mathrm{Im}(P)=\mathrm{Im}(H^\top)$. Let $\check H^\top\in\Rs^{n\times (n-r)}$ be a matrix whose columns form an orthonormal basis for the orthogonal complement of $\mathrm{Im}(P)$ in $\Rs^n$. Then, $\check H^\top \check H+H^\top H=I_n$. Therefore,
\al{ J(\tilde f,g^\star )&= \frac{1}{2}\tr(\tilde P(\check H^\top \check H+H^\top H))=  \frac{1}{2} \tr(\tilde P H^\top H)\nn\\
&=\frac{1}{2} \Es_{\tilde f}[\|H(x-g(y))\|^2].\nn}
This means that the maximization problem can be formulated by reducing the dimension of the random vector $z$, i.e. we take
\al{z_{\mathfrak{r}}:= \left[\begin{array}{c}{   x_{\mathfrak{r} }} \\  y \end{array}\right]=U_{\mathfrak{r}}z, \quad U_{\mathfrak{r}}:=
\left[\begin{array}{cc}
H & 0 \\
0 & I
\end{array}\right],\nn}
and $z=U_{\mathfrak{r}}^\top z_{\mathfrak{r}}$. The nominal and actual reduced densities are, respectively,
$f_{\mathfrak{r}}\sim \mathcal N( m_z^{\mathfrak{r}}, K_z^{\mathfrak{r}})$ and
$\tilde{f}_{\mathfrak{r}} \sim \mathcal N(\tilde m_z^{\mathfrak{r}},\tilde K_z^{\mathfrak{r}})$ where $m_z^{\mathfrak{r}}= U_{\mathfrak{r}} m_z$, $\tilde m_z^{\mathfrak{r}}= U_{\mathfrak{r}}\tilde m_z$,
$K_z^{\mathfrak{r}}= U_{\mathfrak{r}} K_zU_{\mathfrak{r}}^\top$,
$\tilde K_z^{\mathfrak{r}}= U_{\mathfrak{r}} \tilde K_zU_{\mathfrak{r}}^\top$.
Accordingly, the maximization of $J(\cdot,g^0)$ is equivalent to the maximization  $J(\tilde f_{\mathfrak{r}},g_{\mathfrak{r}}^0)=\Es_{\tilde f_{\mathfrak{r}}}[\| x_{\mathfrak{r}} -g^0_{\mathfrak{r}}(y)\|^2]$ where $\tilde f_{\mathfrak{r}}\in \mathcal B_{\mathfrak{r}}:=\{\tilde{f}_{\mathfrak{r}} \sim \mathcal N(\tilde m_z^{\mathfrak{r}},\tilde K_z^{\mathfrak{r}})\hbox{ s.t. } D(\tilde{f}_{\mathfrak{r}}, f_{\mathfrak{r}}) \leq c\}$, $g^0_{\mathfrak{r}}:=Hg^0$ and
\begin{equation*}
\begin{aligned} D(\tilde f_{\mathfrak{r}}, {f_{\mathfrak{r}}})=& \frac{1}{2}\left[(\Delta m^{\mathfrak{r}}_z)^{\top} (K^{\mathfrak{r}}_{z})^{-1}  \Delta m^{\mathfrak{r}}_z +\tr((K^{\mathfrak{r}}_{z})^{-1} \tilde{K}^{\mathfrak{r}}_{z}) \right.\\ &\left.-\ln \operatorname{det} ((K^{\mathfrak{r}}_{z})^{-1} \tilde{K}^{\mathfrak{r}}_{z})-(r+p)\right]\end{aligned}
\end{equation*}
where $\Delta m^{\mathfrak{r}}_{z}=\tilde m^{\mathfrak{r}}_z-m^{\mathfrak{r}}_z$.
Since $K_z^{\mathfrak{r}}>0$ and  $\tilde K_z^{\mathfrak{r}}>0$, by Theorem 1 in \cite{robustleastsquaresestimation, ROBUST_STATE_SPACE_LEVY_NIKOUKHAH_2013} (see also \cite{yi2020low}), we have that the maximizer $\tilde f^0_{\mathfrak{r}}$ has mean $m_z^{\mathfrak{r}}$ and covariance matrix
\begin{equation} \label{tildeK_z}
\tilde{K}^{\mathfrak{r}}_{z}=\left[\begin{array}{cc}\tilde K_x^{\mathfrak{r}}&  K_{xy}^{\mathfrak{r}} \\  K_{yx}^{\mathfrak{r}} & K_y\end{array}\right] \nn
\end{equation} where $\tilde K_x^{\mathfrak{r}}= \tilde P^{\mathfrak{r}}+\tilde K_{xy}^{\mathfrak{r}} K_y^{-1} \tilde K_{yx}^{\mathfrak{r}}$ with
\begin{equation} \label{tilde_P}
\tilde P^{\mathfrak{r}}=\left( (P^{\mathfrak{r}})^{-1}-\lambda^{-1} I_r\right)^{-1},\nn
\end{equation} $P^{\mathfrak{r}}:= HPH^\top >0$  and $\lambda >\sigma_{max}( P^{\mathfrak{r}})>0$ is the unique solution to
\al{\label{eq_lambda}\frac{1}{2}\{ \ln \det P^{\mathfrak{r}} -\ln \det \tilde P^{\mathfrak{r}} +\tr((P^{\mathfrak{r}})^{-1}\tilde P^{\mathfrak{r}}-I_r)\}=c. } We conclude that the least favorable density $\tilde f^0(z)$ has mean and covariance matrix as in (\ref{def_m_K_tilde}). Moreover,
\al{\tilde K_x=H^\top \tilde K_x^{\mathfrak{r}}H = \underbrace{H^\top \tilde P^{\mathfrak{r}}H }_{=:\tilde P}  +\underbrace{ H^\top  K_{xy}^{\mathfrak{r}}}_{=K_{xy}}K_y^{-1}  K_{yx}^{\mathfrak{r}} H \nn}
and \al{\tilde P=H^\top((P^{\mathfrak{r}})^{-1}-\lambda^{-1}I_r)^{-1}H=(P^+-\lambda^{-1}H^\top H)^+.\nn}
Finally, Equation (\ref{eq_lambda}) can be written in terms of $P$ and $\tilde P$:
\al{ \frac{1}{2}\{ \ln {\det}^{+} P  -\ln {\det}^{+} \tilde P  +\tr(P ^{+}\tilde P -I_r)\}=c\nn}
where $\lambda >\sigma_{max}(H P H^\top )=\sigma_{max}(P   ).$

\qed \\

\begin{remark} If $P>0$ then $f$ is a non-degenerate density. In such a case, Theorem \ref{theo1} still holds and: the pseudo inverse is replaced by the inverse; moreover, $H^{ \top} H$ becomes the identity matrix. Therefore, we recover the robust static estimation problem proposed in \cite{robustleastsquaresestimation}.
\end{remark}

\begin{remark}
 In Problem (\ref{robust_p}) we could add the constraint $\rank(G)\leq q$, with $q<r$, i.e. the estimator $g(y)$ represents a layered linear feedforward neural network with one low-dimensional hidden layer. Following the reasonings in \cite{baldi1989neural}, it is possible to characterize $g^0$ in terms of $\tilde f$, however, it is not possible to characterize $\tilde f^0$ in terms of $g^0$. In plain words, the introduction of such a constraint does not allow to characterize the saddle point for Problem  (\ref{robust_p}).
Most importantly, this constraint destroys the convexity of the set $\mathcal G$, so that the saddle point approach to the minimax problem is not applicable.
 \end{remark}

\section{Low-rank robust Kalman filter}\label{sec_3}
We consider a nominal Gauss-Markov state space model of the form
\begin{align}\label{s_space}
\left\{ \begin{array}{rl} x_{t+1}=&A_t x_t+B_t v_t \\ y_t=& C_t x_t+D_t v_t\end{array}\right.
\end{align}
where $x_t\in \Rs^{n}$ is the state vector, $y_t \in \Rs^{p}$ is the measurement vector, $A_t\in\Rs^{n\times n}$, $B_t\in\Rs^{n\times (m+p)}$, $C_t\in\Rs^{p\times n}$, $D_t\in\Rs^{p\times (m+p)}$, $v_t\in \Rs^{m+p}$ is normalized  WGN and $r_t+p:=\rank([\, B_t^{\top}\; D_t^{\top}\, ]^{\top})\leq n+p$. Without loss of generality we can assume that $[\, B_t^{\top}\; D_t^{\top}\,]^{\top}\in \Rs^{(n+p)\times (r_t+p)}$ and it is a full column rank matrix and $v_t\in\mathbb R^{r_t+p}$. Indeed, if this is not the case there always exists a lower dimensional equivalent representation of $v_t$ for which the previous hypothesis holds. Finally, we make the standard assumption that $D_tD_t^{\top}$ is positive definite. We define $z_t:=[\, x_{t+1}^{\top}\; y_t^{\top}\,]^{\top}$. Let $f_0(x_0)$ and $\tilde{\phi}_t(z_t|x_t)$ denote the probability density of $x_0$ and $z_t$ conditioned on $x_t$, respectively. An equivalent representation of model (\ref{s_space}) in the finite interval $0\leq t\leq T$ is given by the nominal joint probability density
 \begin{equation*}
f\left(X_{T+1}, {Y}_{T}\right)=f_{0}\left(x_{0}\right) \prod_{t=0}^{T} \phi_{t}\left(z_{t} | x_{t}\right)
\end{equation*}
where
\begin{equation*}
X_{T+1}=\left[\begin{array}{c}{x_{0}} \\ {\vdots} \\ {x_{t}} \\ {\vdots} \\ {x_{T+1}}\end{array}\right] \quad Y_{T}=\left[\begin{array}{c}{y_{0}} \\ {\vdots} \\ {y_{t}} \\ {\vdots} \\ {y_{T}}\end{array}\right].
\end{equation*}
It is worth noting that $\phi_t(z_t|x_t)$ is a degenerate conditional probability density in the case that $r_t < n$. Let $\mathcal A_t$ denote the support of $\phi_t$. Notice that $\mathcal A_t$ depends on $x_t$. We assume that the actual model over the interval $0\leq t\leq T$
is described by the joint probability
\begin{equation*}
\tilde{f}\left(X_{T+1}, {Y}_{T}\right)=\tilde{f}_{0}\left(x_{0}\right) \prod_{t=0}^{T} \tilde{\phi}_{t}\left(z_{t} | x_{t}\right).
\end{equation*} Moreover, we assume that the support of $\tilde \phi_t(z_t|x_t)$ is $\mathcal A_t$. In this way $f(X_{T+1},Y_T)$ and $\tilde  f(X_{T+1},Y_T)$ have the same support, say $\mathcal A$. Accordingly, we can measure their mismatch through the KL divergence
\begin{equation*}
D(\tilde{f}, f) =\int_{\mathcal A} \tilde{f} \ln\left( \frac{\tilde{f}}{f} \right)d X_{T+1}dY_T.
\end{equation*}
Following reasonings similar to the ones in \cite{ROBUST_STATE_SPACE_LEVY_NIKOUKHAH_2013}, it is not difficult to see that
\begin{equation}\label{D_KL_decomp}
  D(\tilde{f}, f)
 =D(\tilde f_0,f_0)+ \sum_{t=0}^{T} D(\tilde \phi_t,\phi_t)
\end{equation}
where
\begin{align*}
D(&\tilde \phi_t,\phi_t)= \tilde \Es\left[\ln\left(\frac{\tilde \phi_t}{\phi_t}\right)\right]\nn\\
& =\int_{ \tilde{ \mathcal A}_t} \int_{\mathcal A_t} \tilde f_t(x_t) \tilde \phi_t(z_t|x_t) \ln\left(\frac{\tilde \phi_t}{\phi_t}\right)d z_t d x_t,
\end{align*}
and $\tilde f_t(x_t)$ denotes the actual marginal density of $x_t$ whose support is denoted by $\tilde{\mathcal A}_t$.

Assume to know the nominal model (\ref{s_space}), and thus $f(X_{T+1},Y_T)$, while the actual one does not. In view of the decomposition in
(\ref{D_KL_decomp}) and given $Y_{t-1}$, we can assume that $\tilde \phi_t$ belongs to the following convex ambiguity set:
\begin{equation*}
\mathcal B_t:=\left\{\, \tilde{\phi}_{t} \hbox{ s.t. }\tilde{\Es}\left[\ln \left(\frac{\tilde{\phi}_{t}\left(z_{t} | x_{t}\right)}{\phi_{t}\left(z_{t} | x_{t}\right)}\right) \bigg| Y_{t-1}\right] \leq c_{t}\right\}
\end{equation*}
where
\begin{equation*}
\begin{aligned}&\tilde{\Es}\left[\ln \left(\frac{\tilde{\phi}_{t}\left(z_{t} | x_{t}\right)}{\phi_{t}\left(z_{t} | x_{t}\right)}\right) \bigg| Y_{t-1}\right] \\ &:= \int_{\check{\mathcal A}_t}\int_{\mathcal A_t} \tilde{\phi}_{t}\left(z_{t} | x_{t}\right) \tilde{f}_{t}\left(x_{t} | Y_{t-1}\right) \ln \left(\frac{\tilde{\phi}_{t}\left(z_{t} | x_{t}\right)}{\phi_{t}\left(z_{t} | x_{t}\right)}\right) d z_{t} d x_{t}
 \end{aligned}
\end{equation*} and $\check{\mathcal A}_t$ denotes the support of $\tilde f_t(x_t|Y_{t-1})$. In plain words, the ambiguity set is expressed incrementally by forming a ball, in the KL topology, about the nominal model and placing an upper bound $c_t$ on its radius. Then, we consider the robust estimator of $x_{t+1}$ given $Y_t$  solving the following dynamic minimax game:
\begin{equation}
 \label{minimax} \hat x_{t+1} =\underset{g_t \in \mathcal{G}_{t}}{\mathrm{argmin}}\max_{\tilde{\phi}_{t} \in \mathcal{B}_{t}} J_t(\tilde {\phi}_t,g_t)
\end{equation}
where \begin{equation*}
 \begin{aligned}
J_t(\tilde {\phi}_t,g_t)=\frac{1}{2}\tilde{\mathbb{E}}\left[\left\|   x_{t+1}-g_t\left(y_{t}\right) \right\|^{2} | \mathrm{Y}_{t-1}\right]\\
=\frac{1}{2} \int_{\check A_t}\int_{\mathcal A_t}\left\|  x_{t+1}-g_{t}\left(y_{t}\right) \right\|^{2} \tilde{\phi}_{t}(z_{t} | x_{t})  \\
 \times \tilde{f}_{t}\left(x_{t} | Y_{t-1}\right) d z_{t} d x_{t};
\end{aligned}
\end{equation*}$\mathcal{G}_{t}$ denotes the class of estimators with finite second-order moments with respect to all the densities $\tilde{\phi}_{t}\left(z_{t} | x_{t}\right) \tilde{f}_t(x_t|Y_{t-1})$ such that $\tilde{\phi}_{t} \in \mathcal{B}_{t}$.  Note that $\tilde \phi_t$ must satisfy the constraint:
\begin{equation}
I_{t} (\tilde{\phi}_{t} ) \triangleq \int_{ \check A_t}\int_{\mathcal A_t} \tilde{\phi}_{t}\left(z_{t} | x_{t}\right) \tilde{f}_{t}\left(x_{t} | Y_{t-1}\right) d z_{t} d x_{t}=1. \label{I}
\end{equation}
\begin{lemma}
\label{lemma2}
For a fixed estimator $g_t \in \mathcal{G}_t$, the density $\tilde {\phi}_{t}\left(z_{t} | x_{t}\right) \in \mathcal{B}_{t}$ {having support $\mathcal A_t$} that maximizes the objective function
$$  J_t(\tilde {\phi}_t,g_t)=\tilde{\mathbb{E}}\left[\left\|   x_{t+1}-g_t  (y_{t}  )\right\|^{2} | \mathrm{Y}_{t-1}\right]$$
under the constraint $D(\tilde{\phi}_{t}, \phi_{t}) \leq c_{t}$ is given by
\begin{equation} \label{phi_0}
\tilde{\phi}_{t}^{0}=\frac{1}{M_{t}\left(\lambda_{t}\right)} \exp \left(\frac{1}{2 \lambda_{t}}\left\|  x_{t+1}-g_{t} (y_{t})\ \right\|^{2}\right) \phi_{t}.
\end{equation} Moreover, $M_t(\lambda_t) $ is the normalizing constant such that (\ref{I}) holds. Finally, {for $c_t>0$ sufficiently small, there exists} a unique $\lambda_t>0$ such that $D(\tilde \phi^0_t, \phi_t)=c_t $.
\end{lemma}
\IEEEproof
For a given $g_t$, the Lagrangian for the constrained optimization problem takes the form:
\begin{equation*}
\begin{aligned} L_{t}& (\tilde{\phi}_{t}, \lambda_{t}, \mu_{t} )\\
=&  J_{t} (\tilde{\phi}_{t}, g_t ) +\lambda_{t} (c_{t}-D_{t}(\tilde{\phi}_{t}, \phi_{t} ) ) +\mu_{t} (1-I_{t}(\tilde{\phi}_{t} ) )\\
=&\frac{1}{2} \int_{\check{\mathcal A}_t}\int_{\mathcal A_t}\left\|    x_{t+1}-g_{t} (y_{t} ) \right\|^{2} \tilde{\phi}_{t} \tilde{f}_{t}\left(x_{t} | Y_{t-1}\right) d z_{t} d x_{t}\\
&+\lambda_t\left(c_t-\int_{\check{\mathcal A}_t}\int_{\mathcal A_t} \tilde{\phi}_{t}\tilde{f}_{t}\left(x_{t} | Y_{t-1}\right) \ln \left(\frac{\tilde{\phi}_{t}}{\phi_{t}}\right) d z_{t} d x_{t}\right)
\\&+\mu_{t} \left(1-\int_{\check{\mathcal A}_t}\int_{\mathcal A_t} \tilde{\phi}_{t} \tilde{f}_{t}\left(x_{t} | Y_{t-1}\right) d z_{t} d x_{t} \right)
  \end{aligned}
\end{equation*}
where $\lambda_t\geq 0$. In order to maximize $ L_{t}(\tilde{\phi}_{t}, \lambda_{t}, \mu_{t})$ with respect to $\tilde \phi _t$, we need to prove {it is concave} in $\tilde \phi_t$. The first and second variation of $ L_{t}$ with respect to $\tilde \phi_t$ along the direction $\delta \tilde{\phi}_{t} $ are, respectively,
\begin{equation}
\begin{aligned} \delta L_t&(\tilde{\phi}_{t}, \lambda_{t}, \mu_{t} ;\delta \tilde{\phi}_{t})\\
=&{\frac{1}{2} \int_{\check{\mathcal A}_t} \int_{\mathcal A_t} \left\|  x_{t+1}-g_t\left(y_{t} \right)\right\|^{2} \delta \tilde{\phi}_{t}\tilde{f}_{t}\left(x_{t} | Y_{t-1}\right) d z_t d x_t}\\
 &+\lambda_t\left(-\int_{\check{\mathcal A}_t} \int_{\mathcal A_t} \delta \tilde{\phi}_{t} \tilde{f}_{t} \left(x_{t} | Y_{t-1}\right)\ln\left( \frac{\tilde{\phi}_{t}}{\phi_{t}}\right) d z_t d x_t  \right.\\
 &\left.  -\int_{\check{\mathcal A}_t}\int_{\mathcal A_t}  \tilde{f}_{t}\left(x_{t} | Y_{t-1}\right)  \delta \tilde{\phi}_{t} d z_t d x_t\right)\\
& -\mu_t \int_{\check{\mathcal A}_t}\int_{\mathcal A_t} \delta \tilde{\phi}_{t} \tilde{f}_t\left(x_{t} | Y_{t-1}\right)  d z_t d x_t\\
 \delta L&^{2}(\tilde{\phi}_{t}, \lambda_{t}, \mu_{t} ; \delta \tilde{\phi}_{t}^{2})\nn\\ &=-\lambda_t\int_{\check{\mathcal A}_t}\int_{\mathcal A_t} \delta \tilde{\phi}^2_{t} \tilde{f}_{t}\left(x_{t} | Y_{t-1}\right) \frac{1}{\tilde{\phi}_{t}}   d z_t d x_t.
\end{aligned}
\end{equation}
Notice that $\delta \tilde{\phi}_{t}(\cdot |x_t)$ is a function whose support is $\mathcal A_t$. Then, under the assumption that $\phi_t>0$ in $\mathcal A_t$ and $\lambda_t>0$, we have that the second variation is negative for any $\delta \tilde{\phi}_{t} \neq 0$. So $L_t$ is strictly concave in $\tilde{\phi}_{t}$. This implies that the point of maximum is given by imposing the stationarity condition  $\delta L(\tilde{\phi}_{t}, \lambda_{t}, \mu_{t} ; \delta \tilde{\phi}_{t})=0$
 for all functions $\delta \tilde{\phi}_t$ (having support $\mathcal A_t$) which implies
\begin{equation*}
\begin{aligned}
\left[\frac{1}{2} \left\|  x_{t+1}-g_t\left(y_{t}\right) \right\|^{2}-\lambda_t \ln \left(\frac{\tilde{\phi}_{t}}{\phi_{t}}\right)-\lambda_t-\mu_t\right]\\
\times  \tilde{f}_t\left(x_{t} | Y_{t-1}\right)=0.&
\end{aligned}
\end{equation*}
Accordingly, for any $x_t\in \check{\mathcal A}_t$, the maximum $\phi _t^0$ is such that
\begin{equation*}
\ln \left(\frac{\tilde{\phi}_{t}^0}{\phi_{t}}\right)=\frac{1}{2 \lambda_{t}} \left\|  x_{t+1}-g_{t}(y_{t})\right\|^{2}-\ln M_{t}(\lambda _t)
\end{equation*}
where $$
\ln M_{t}(\lambda_t) \triangleq 1+\frac{\mu_{t}}{\lambda_{t}}.
$$
Thus, \begin{equation*}
\tilde{\phi}_{t}^{0}=\frac{1}{M_{t}\left(\lambda_{t}\right)} \exp \left(\frac{1}{2 \lambda_{t}} \left\|  x_{t+1}-g_{t}\left(y_{t} \right)\right\|^{2}\right) \phi_{t}
\end{equation*}
{and $M_t(\lambda_t)$} must satisfy condition (\ref{I}):
\begin{equation*}
\begin{aligned}
\int_{\check{\mathcal A}_t} \int_{\mathcal A_t} \frac{1}{M_t(\lambda_t)} \exp \left(\frac{1}{2 \lambda_t} \left\|  x_{t+1}-g_t\left(y_{t}\right) \right\|^{2}\right)\\
 \times{\phi}_{t} \tilde{f}_{t}\left(x_{t} | Y_{t-1}\right)  d z_t d x_t=1&
\end{aligned}
\end{equation*}
so that
\begin{equation*}
\begin{aligned} M_{t}\left(\lambda_{t}\right)=&\int_{\check{\mathcal A}_t}\int_{\mathcal A_t} \exp \left(\frac{1}{2 \lambda_{t}}\left\|   x_{t+1}-g_{t}\left(y_{t} \right) \right\|^{2}\right) \phi_{t} \\ & \times \tilde{f}_{t}\left(x_{t} | Y_{t-1}\right) d z_{t} d x_{t}. \end{aligned}
\end{equation*}
In order to find the unique Lagrange multiplier $\lambda_t>0$ satisfying the inequality constraint, we consider
\begin{equation*}
\begin{aligned} D&(\tilde{\phi}_{t}^{0}, \phi_{t}) =\int_{\check{\mathcal A}_t} \int_{\mathcal A_t} \tilde \phi^0_t\tilde{f}_{t}\left(x_{t} | Y_{t-1}\right)  \ln \left(\frac{\tilde{\phi}_{t}^{0}}{\phi_{t}} \right)d z_t d x_t \\ & =\int_{\check{\mathcal A}_t}\int_{\mathcal A_t} \left(\frac{1}{2 \lambda_{t}} \left\|  x_{t+1}-g_t\left(y_{t}\right) \right\|^{2}-\ln M_{t}\left(\lambda_{t}\right)\right) \\
&\times \tilde{\phi}_{t}^{0} \tilde{f}_{t} \left(x_{t} | Y_{t-1}\right) d z_t d x_t \\ &=\frac{1}{ \lambda_{t}} J_{t}\left(\tilde{\phi}_{t}^{0}, g_{t}\right)-\ln M_{t}\left(\lambda_{t}\right). \end{aligned}
\end{equation*}
Notice that $$
\begin{aligned}  &\frac{d}{d \lambda_{t}}\ln M_{t}\left(\lambda_{t}\right)\\
 =& \frac{-1}{2 \lambda_{t}^{2}}\frac{1}{M_{t}\left(\lambda_{t}\right)}\int_{\check{\mathcal A}_t} \int_{\mathcal A_t} \exp \left(\frac{1}{2 \lambda_{t}} \left\|  x_{t+1}-g_t\left(y_{t}\right)\right\|^{2}\right)  \\ & \times \phi_{t} \tilde{f}_{t} \left(x_{t} | Y_{t-1}\right) \left\| x_{t+1}-g_t\left(y_{t}\right) \right\|^{2} d z_t d x_t \\  =&  \frac{-1}{2 \lambda_{t}^{2}} \int_{\check{\mathcal A}_t}\int_{\mathcal A_t} \left\| x_{t+1}-g_t\left(y_{t}\right)\right\|^{2}\tilde{f}_{t}\left(x_{t} | Y_{t-1}\right) \tilde  \phi_{t}   d z_t d x_t \\=&- \frac{1}{ \lambda_{t}^{2}} J_t\left(\tilde{\phi}_{t}^{0}, g_t(y_t)\right) \end{aligned}
$$
so that
\begin{equation*}
\kappa_t\left(\lambda_{t}\right) \triangleq D\left(\tilde{\phi}_{t}^{0}, \phi_{t}\right)=-\lambda_{t} \frac{d}{d \lambda_{t}} \ln M_{t}\left(\lambda_{t}\right)-\ln \left(M_{t}\left(\lambda_{t}\right)\right).
\end{equation*}

Then, the derivative of $ \kappa_t\left(\lambda_{t}\right)$ is given by:
\begin{equation*}
\begin{aligned} \frac{d {\kappa}_{t}}{d \lambda_{t}}=&-\lambda_{t} \frac{d^{2}}{d \lambda_{t}^{2}} \ln M_{t}-2 \frac{d}{d \lambda_{t}} \ln M_{t} \\=&-\frac{1}{\lambda_{t}}\left[\frac{d}{d \lambda_{t}}\left(\lambda_{t}^{2} \frac{d}{d \lambda_{t}} \ln M_{t}\right)\right]=\frac{1}{\lambda_{t}} \frac{d}{d \lambda_{t}} J_{t}\left(\tilde{\phi}_{t}^{0}, g_{t}\right) \\=&{ -\frac{1}{4 \lambda_{t}^{3}} \tilde{\Es}\left[\left( \| x_{t+1}-g_{t} (y_{t})\ \|^{2}\right.\right.}\\ & \left.\left.-\tilde{\Es}\left[ \left\| x_{t+1}-g_{t}(y_{t})\right\|^{2} | Y_{t-1}\right]\right)^{2} | Y_{t-1}\right] <0.\end{aligned}
\end{equation*}
Therefore,  ${\kappa_t}\left(\lambda_{t}\right)$ is a monotone decreasing function. Notice that when $\lambda_{t}  \rightarrow \infty$, we have $\tilde{\phi}_{t}^{0} \rightarrow \phi_{t}$, so that ${\kappa_t}(\infty)=0$. Accordingly, if $c_t$ is sufficiently small, then $c_t$ is in the range of ${ \kappa_t}$. Therefore, there exists a unique $\lambda_t>0$ such that ${\kappa_t}(\lambda_t)=c_t$.
\qed \\

Once we get the function $\tilde \phi^0_t$, the estimator $g_t \in \mathcal G_t$ minimizing the objective function $J_t(\tilde \phi_t^0,g_t)$ is given by
\begin{equation}\label{hatx}
\begin{aligned} \hat{x}_{t+1}=g_{t}^{0}\left(y_{t}\right) &=\tilde{\Es}\left[x_{t+1} | Y_{t}\right] \\ &=\int_{ {\mathcal A}_{t}^\star} x_{t+1} \tilde{f}_{t+1}\left(x_{t+1} | Y_{t}\right) d x_{t+1} \end{aligned}
\end{equation}
where \begin{equation*}
\tilde{f}_{t+1}\left(x_{t+1} | Y_{t}\right)=\frac{\int_{\check{\mathcal A_t}} \tilde{\phi}_{t}^{0}\left(z_{t} | x_{t}\right)  \tilde{f}_{t}\left(x_{t} | Y_{t-1}\right) d x_{t}}{\int_{\check{\mathcal A}_t} \int_{\mathcal A_{t}^{\star\star}} \tilde{\phi}_{t}^{0}\left(z_{t} | x_{t}\right)  \tilde{f}_{t}\left(x_{t} | Y_{t-1}\right) d x_{t+1} d x_{t}},
\end{equation*}
$\mathcal A^\star_t$ is the support of $\tilde f_t(x_{t+1}|Y_t)$ and $\mathcal A^{\star\star}_t$ is defined as follows: $x_{t+1}\in \mathcal A^{\star\star}_t$ if and only if there exists at least one $y_t$ for which $ [\,x_{t+1}^{\top}\; y_t^{ \top}\,]^{ \top}\in\mathcal A_t $.

\section{Low-rank robust Kalman filter as a static game} \label{sec_4}
In view of (\ref{hatx}), the optimal estimator $\hat x_{t+1}$ depends on the least favorable density $\tilde{\phi}_{t}^{0}\left(z_{t} | x_{t}\right)$ while from  (\ref{phi_0}), the least favorable density $\tilde{\phi}_{t}^{0}\left(z_{t} | x_{t}\right)$ depends on the optimal estimator $\hat x_{t+1}$. In order to break this endless loop, we assume at time $t$ the $a~ priori $ conditional density of $x_t$ given $ Y_{t-1}$ is
\begin{align}\label{hp_ft}
\tilde{f}_{t}\left(x_{t} | Y_{t-1}\right) \sim \mathcal{N} (\hat{x}_{t}, \tilde P_{t})
\end{align}
and let $\rank(\tilde P_t)=r_t$. In addition, we know the nominal conditional transition probability density function of model (\ref{s_space}) is $\phi_t(z_t|x_t) \sim \mathcal{N}\left(m_{z_t|x_t}, K_{z_t|x_t}\right)$ with$$
m_{z_t|x_t}=\left[\begin{array}{c}{ A_t} \\ { C_t}\end{array}\right] x_{t}, ~~~K_{z_t|x_t}=\left[\begin{array}{cc}{ B_t  B_t^{\top}} & {B_t D_t^{\top}} \\ { D_tB_t ^{\top}} & { D_t  D_t^{\top}}\end{array}\right].$$
As we already noticed, the latter could be degenerate and thus
 $K_{z_t|x_t}$ could be singular. We define the Gaussian pseudo nominal density
\begin{equation*}
{ {f}_{t}\left(z_{t} | Y_{t-1}\right)}=\int_{\check{\mathcal A}_t} \phi_{t}\left(z_{t} | x_{t}\right){ \tilde{f}_{t}\left(x_{t} | Y_{t-1}\right)} d x_{t}.
\end{equation*} Since both $\phi_t(z_t|x_t)$ and $\tilde{f}_{t}\left(x_{t} | Y_{t-1}\right)$ are Gaussian, we obtain ${ {f}_{t}\left(z_{t} | Y_{t-1}\right)} \sim \mathcal{N}\left(m_{z_t}, K_{z_{t}}\right)$ with
\begin{equation*}
m_{z_t}=\left[\begin{array}{l}
A_{t} \\
C_{t}
\end{array}\right] \hat{x}_{t}, \quad {K}_{z_{t}}=\left[\begin{array}{cc}
{K}_{x_{t+1}} & K_{x_{t+1} y_{t}} \\
K_{y_{t} x_{t+1}} & K_{y_{t}}
\end{array}\right].
\end{equation*}
where the parametric form of $K_{z_t}$ is given by
\begin{equation}\label{K_z}
K_{z_{t}}=\left[\begin{array}{c}
A_{t} \\
C_{t}
\end{array}\right] \tilde P_{t}\left[\begin{array}{cc}
A_{t}^{ \top} & C_{t}^{ \top}
\end{array}\right]+\left[\begin{array}{c}
B_{t} \\
D_{t}
\end{array}\right]\left[\begin{array}{cc}
B_{t}^{ \top} & D_{t}^{\top}
\end{array}\right].
\end{equation}
Then, we have $$ K_{y_t}= C_t\tilde P_t C_t^{\top}+D_tD_t^{\top}\geq D_tD_t^{\top}>0.$$  Next, we define the actual density
\begin{align*}
\tilde{f}_{t}\left(z_{t} | Y_{t-1}\right)=\int_{{\check{\mathcal A}}_t}\tilde{\phi}_{t}\left(z_{t} | x_{t}\right) {\tilde{f}_{t}\left(x_{t} | Y_{t-1}\right)} d x_{t}.
\end{align*}
which has the same support of  $f_t(z_t|Y_{t-1})$, which is denoted by $\breve{\mathcal A}_t$. Moreover, on the basis of Equation (\ref{phi_0}), it is easy to show that if $\phi_t(z_t|x_t)$ is Gaussian, the least favorable density
\begin{equation*}
\begin{aligned}&\tilde{f}_{t}^{0}\left(z_{t} | Y_{t-1}\right)\\
&{=\frac{1}{M\left(\lambda_{t}\right)} \exp \left(\frac{1}{2 \lambda_{t}} \left\| x_{t+1}-g_{t}\left(y_{t}\right) \right\|^{2}\right) {f}_{t}\left(z_{t} | Y_{t-1}\right)}
\end{aligned}\end{equation*}
is also Gaussian. Furthermore, it is not difficult to see that  $D(\tilde \phi_t^0, \phi_t)= D( \tilde f_t^0,  f_t)$. Hence, the least favorable density of $z_t$ conditioned on $Y_{t-1}$ belongs to the following ambiguity set
\begin{align*}
\mathcal{\bar B}_t=\{ \tilde f_t\sim \mathcal N(\tilde m_{z_t},\tilde K_{z_t}) ~~  s.t.\;\;   D(\tilde f_t,  f_t) \leq c_t\}.
\end{align*}
Accordingly, the dynamic minimax game in (\ref{minimax}) is equivalent to the following static minimax game
\begin{equation} \label{minimax_dyn}
\hat x_{t+1} =\underset{g_{t} \in \mathcal{G}_{t}}{\mathrm{argmin}} \max _{\tilde{f}_{t} \in \bar{\mathcal{B}}_{t}} \bar J_{t}(\tilde f_{t}, g_{t} )
\end{equation}
where the least favorable model is now characterized by the conditional density $\tilde f_t(z_t|Y_{t-1})$ and $$
\bar{J}_{t}(\tilde{f}_{t}, g_{t})=\int_{\breve{\mathcal A}_t}{ \left\| x_{t+1}-g_{t}\left(y_{t}\right) \right\|^{2}} \tilde{f}_{t}\left(z_{t} | Y_{t-1}\right) \mathrm{d} z_{t}.
$$  It is worth noting that the minimizer of $\bar J_t( f_t,\cdot )$ is the Bayesian estimator and the state prediction error $x_{t+1}-g_t(y_t)$ is Gaussian distributed with zero mean and covariance matrix
\begin{align}\label{Pt1}
P_{t+1} &:=K_{x_{t+1}}-K_{x_{t+1}, y_{t}} K_{y_{t}}^{-1} K_{y_{t}, x_{t+1}}
\end{align}
which could be singular.
In such a case,  \al{\label{cond_pred_err}H_t^\top H_t(x_{t+1}-g_t(y_t))=x_{t+1}-g_t(y_t)} where $H_t^\top$ is a matrix
whose columns form an orthonormal basis for $\mathrm {Im} (P_{t+1})$. On the other hand, the prediction error under the least favorable model is Gaussian with zero mean and covariance matrix $\tilde P_{t+1} $ such that $\mathrm{Im}(\tilde P_{t+1})=\mathrm{Im}(P_{t+1})$ because $\mathrm{Im}(\tilde K_{z_t})= \mathrm{Im}(K_{z_t})$. Therefore, condition (\ref{cond_pred_err}) still holds and the least favorable probability density of $z_t$ given $x_t$ can be written as
\begin{equation} \label{phi_0new}
\tilde{\phi}_{t}^{0}=\frac{1}{M_{t}\left(\lambda_{t}\right)} \exp \left(\frac{1}{2 \lambda_{t}}\left\| H_t (x_{t+1}-g_{t}\left(y_{t})\right)\right\|^{2}\right) \phi_{t}.
\end{equation}

\begin{remark}  Condition (\ref{hp_ft}) means that, using the terminology coined by Hansen and Sargent \cite{ROBUSTNESS_HANSENSARGENT_2008}, the maximizer in (\ref{minimax_dyn}) is operating under commitment, i.e. the maximizer is required to commit all the least favorable model components at early stages  with the estimating player.
\end{remark}

At this point the solution to (\ref{minimax_dyn}) is given by Theorem  \ref{theo1}: substituting $f$, $\tilde f$, $g$, $H$ with $f_t$, $\tilde f_t$, $g_t$, $H_t$, respectively, it is not difficult to see that all the assumptions are satisfied. Then, the least favorable density is $\tilde{f}_{t}^{0}(z_{t} | Y_{t-1}) \sim \mathcal{N}(m_{z_t}, \tilde{K}_{z_t})$ where $$
\tilde{K}_{z_{t}}=\left[\begin{array}{cc}
\tilde{K}_{x_{t+1}} & K_{x_{t+1}, y_{t}} \\
K_{y_{t}, x_{t+1}} & K_{y_{t}}
\end{array}\right].
$$
The nominal posterior covariance of $x_{t+1}$ given $Y_{t}$ has been defined in (\ref{Pt1}).
Accordingly, the least favorable posterior covariance of $x_{t+1}$ given $Y_{t}$  is
\begin{align*}
\tilde P_{t+1} &=\tilde{K}_{x_{t+1}}-K_{x_{t+1}, y_{t}} K_{y_{t}}^{-1} K_{y_{t}, x_{t+1}}.
\end{align*}
Moreover, in view of Equation (\ref{tilde_P}), we have \begin{equation}\label{def_P_tilde}
\tilde P_{t+1}=\left(P_{t+1}^{+}-\lambda_{t+1}^{-1} H_t^{\top} H_t\right)^{+}
\end{equation}
where $\lambda_t>\sigma_{max}( P_{t+1} )$ and the Lagrange multiplier $\lambda_t$ is selected by the following equation
\begin{equation}\label{def_gamma}
 \begin{aligned}
  \gamma(P_{t+1},\lambda_t):=&\frac{1}{2}\left\{\operatorname{lndet}^+(P_{t+1})-\operatorname{lndet}^+(\tilde P_{t+1})\right.\\
&\left.+\tr\left[P_{t+1}^+\tilde P_{t+1}-I_{r_{t+1}}\right]\right\}=c_t.
\end{aligned}
\end{equation}
Furthermore, the robust estimator is
\begin{equation*}
\hat{x}_{t+1}=A_{t} \hat{x}_{t}+G_{t}\left(y_{t}-C_{t} \hat{x}_{t}\right).
\end{equation*}
where $G_t=K_{x_{t+1},y_t}K^{-1}_{y_t}$. Then, in view of Equation (\ref{K_z}), we know $G_t$ and $P_{t+1}$ take the parametric form
\begin{equation*}
G_{t}= (A_{t} \tilde P_{t} C_{t}^{ \top}+B_{t} D_{t}^{\top})(C_{t} \tilde P_{t} C_{t}^{\top}+D_{t} D_{t}^{\top})^{-1}
\end{equation*}
\begin{equation*}
P_{t+1}=A_{t} \tilde P_{t} A_{t}^{\top}-G_{t} (C_{t} \tilde P_{t} C_{t}^{\top}+D_{t} D_{t}^{\top}) G_{t}^{\top}+B_{t} B_{t}^{\top}.
\end{equation*}

\begin{algorithm}[h]
  \caption{Low-rank robust Kalman filter at time $t$}
  \begin{algorithmic}[1]
    \Require
      $y_t$, $\hat x_t$, $\tilde P_t$, $c_t$
      \State $ G_{t}=(A_{t} \tilde P_{t} C_{t} ^{\top}+B_{t} D_{t}^{\top})(C_{t} \tilde P_{t} C_{t}^{ \top}+D_{t} D_{t}^{\top})^{-1}$
      \State $\hat{x}_{t+1}=A_{t} \hat{x}_{t}+G_{t}\left(y_{t}-C_{t} \hat{x}_{t}\right)$
       \State \label{RIcc_step}$P_{t+1}=A_{t} \tilde P_{t} A_{t}^{\top}-G_{t}\left(C_{t} {\tilde P_{t} }C_{t}^{ \top}+D_{t} D_{t}^{ \top}\right) G_{t}^{\top}+B_{t} B_{t}^{ \top}$
      \State Select $H_t^{\top}H_t$ as projection matrix with image $ \mathrm {Im} (P_{t+1})$
      \State \label{RIcc_stepfind} Find {$\lambda_t^{-1}$} s.t. $ \gamma( P_{t+1},\lambda_t)=c_t$
      \State \label{RIcc_step2} $\tilde P_{t+1}=(P_{t+1}^{+}-\lambda_{t}^{-1} H_t^{\top} H_t)^{+}$
  \end{algorithmic}\label{code:recentEnd}
\end{algorithm}

Algorithm \ref{code:recentEnd} is the summary of the resulting low-rank robust Kalman filter. Note that steps \ref{RIcc_step}-\ref{RIcc_step2} correspond to a distorted Riccati iteration involving covariance matrices which are possibly singular. In addition, $\theta_t:=\lambda_t^{-1}$ is the time-varying risk sensitivity parameter. It is worth noticing that when $c_t=0$ the actual model coincides with the nominal one, so that $\theta_t=0$, and thus $P_t=\tilde P_t$. Then, step \ref{RIcc_step} in Algorithm \ref{code:recentEnd} becomes the usual Riccati equation which means we obtain the standard Kalman filter.
\begin{remark}  {The computational complexity of Algorithm \ref{code:recentEnd} is now discussed. Steps 1-3 have the same complexity of the standard Kalman filter, that is
$$O(n^3)+O(n(r+p)p)+O(p^2(r+p))+O(n^2(r+p)).$$
Step 4 and Step 6 have complexity $O(n^3)$. In regard to Step 5, the computation of $\lambda_t^{-1}\in (0,\sigma_{max}(P_{t+1})^{-1})$ is accomplished by a bisection method, see Algorithm 2 in \cite{zenere2018coupling}.  The complexity to evaluate $\gamma(P_{t+1},\lambda_t)$ is $O(n^3)$, thus the complexity of Step 5 is $O(\log_2(\frac{\sigma_{max}(P_{t+1})^{-1}}{\epsilon})n^3)$ where $\epsilon>0$ is the selected accuracy, i.e. the solution found satisfies the condition $| \gamma( P_{t+1},\lambda_t)-c_t|\leq \epsilon$. We conclude that the computational complexity of Algorithm \ref{code:recentEnd} is
  \al{O(n^3)+&O(n(r+p)p)+O(p^2(r+p))+O(n^2(r+p))\nn\\ &+O(\log_2(\frac{\sigma_{max}(P_{t+1})^{-1}}{\epsilon})n^3)\nn.}}
  \end{remark}{\em Low-rank risk-sensitive filtering}: The minimax problem in (\ref{minimax_dyn}) can be relaxed by imposing the constraint $\tilde f_t\in \bar {\mathcal B}_t$ through a penalty term:
\begin{equation} \label{minimax_dyn_relax}
\hat x_{t+1} =\underset{g_{t} \in \mathcal{G}_{t}}{\mathrm{argmin}}  \, \bar J_{t}(\tilde f_{t}, g_{t} )-\theta D(\tilde f_t,f_t)
\end{equation}
where $\theta>0$ represents the risk sensitivity parameter and it can be understood as a regularization parameter. It is not difficult to see that the solution to (\ref{minimax_dyn_relax}) is the one in Algorithm \ref{code:recentEnd} where Steps \ref{RIcc_stepfind}-\ref{RIcc_step2} are replaced by
$$ \tilde P_{t+1}=(P_{t+1}^+-\theta H_t^{\top} H_t)^+.$$
Such an estimator represents the extension to the singular case of the well known risk sensitive filtering problem, \cite{boel2002robustness,RISK_WHITTLE_1980}.

\section{The least favorable model} \label{sec_5}
In this section, we characterize the least favorable model corresponding to (\ref{minimax}). Let $e_{t+1}=x_{t+1}-\hat x_{t+1}$ denote the state prediction error. By (\ref{s_space}), we have
\begin{equation}\label{e_t+1}e_{t+1}=\left( A_t- G_{t} C_t\right) e_{t}+\left( B_t- G_{t}  D_t\right) v_{t}.\end{equation}
It is worth noticing that the driving noise $v_t$ is independent from $e_t$ under the nominal model. In view of (\ref{phi_0new}), we obtain the least favorable density  \begin{equation}
\tilde{\phi}_{t}^{0}=\frac{1}{M_{t}\left(\lambda_{t}\right)} \exp \left(\frac{1}{2 \lambda_{t}}\left\| H_t e_t\right\|^{2}\right) \phi_{t}.
\end{equation}
The latter is not a normalized density, which means that the hostile player could have the opportunity to change retroactively the least favorable density of $x_t$.

Accordingly, it is possible to find the least favorable density with respect to the driving noise $v_t$. From model (\ref{s_space}), it is not difficult to see that
\begin{equation}
\left[\begin{array}{l}  B_t \\  D_t \end{array}\right]v_t=z_{t}-\left[\begin{array}{l} A_t \\ C_t \end{array}\right] x_{t} .
\end{equation}
where $z_t$ belongs to the affine space  $$\mathbb{S}=\{ [\,A_t^\top \; C_t^\top \,]^\top x_t+ \bar w_t,~ \bar w_t \in \operatorname{Im}(\Gamma_t)\}$$ and $\Gamma_t:=[\, B_t^\top\; D_t^\top\,]^\top$ is full column rank. Accordingly, there is a one-to-one correspondence between $v_t\in \mathbb R^{r_t+p}$ and $z_t\in \mathbb  S$ through the isomorphism
\al{ &\mathbb R^{r_t+p}  \rightarrow  \mathbb S\nn\\ & v_t \mapsto  (\Gamma_t^{ \top} \Gamma_t)^{-1}\Gamma_t^{ \top} \left( z_{t}-\left[\begin{array}{l} A_t \\ C_t \end{array}\right] x_{t}\right).\nn}

Therefore, instead of characterizing the least favorable density $\tilde{\phi}_{t}^{0}(z_{t} | x_{t}) $ directly, we characterize it in terms of $v_t$. The nominal density is denoted by $\psi_{t}(v_{t})$ and \begin{equation*}\psi_{t}\left(v_{t}\right) \sim \exp \left(-\left\|v_{t}\right\|^{2} / 2\right)\end{equation*}
where $\sim$  means that the two terms are the same up to constant scale factors. Then, we make the guess that least favorable density of $v_t$ depends on the state prediction error. Hence, the least favorable density is denoted by $\tilde \psi_{t}(v_{t} | e_{t})$. We define it over the interval $[t+1,T]$. Then,
\begin{equation*}\begin{aligned}
 \prod_{s=t+1}^{T} &\exp \left(\frac{\left\|{H}_t  e_{s+1} \right\|^{2}}{2\lambda_{s}}\right) \psi_{s}\left(v_{s}\right)\\
&\sim \exp \left(\frac{\left\| e_{t+1}\right\|_{ \Omega_{t+1}^{+}}^{2}}{2}\right) \prod_{s=t+1}^{T} \tilde {\psi}_{s}\left(v_{s} |  e_{s}\right)
\end{aligned}\end{equation*}
where the term $\exp (\left\| e_{t+1}\right\|_{ \Omega_{t+1}^{+}}^{2})$ indicates the cumulative error at time $t+1$. Note that $\Omega_{t+1}^+$ could be singular. Decreasing the time index $t$ by 1 we obtain:
\begin{equation*}\begin{aligned}
\exp &\left(\frac{\left\| H_t  e_{t+1}\right\|^2}{2\lambda_{t}}\right)\psi_{t}\left(v_{t}\right) \\
&\sim\exp \left(\frac{\left\|e_{t}\right\|_{ \Omega_{t}^{+}}^{2}-\left\|e_{t+1}\right\|_{ \Omega_{t+1}^{+}}^{2}}{2}\right) \tilde \psi_{t}\left(v_{t}| e_t\right) .
\end{aligned}\end{equation*}
Hence, we have
\begin{equation*}\begin{aligned}
\tilde \psi_{t}\left(v_{t}|e_t\right) \sim \exp  \left(\frac{\left\| e_{t+1}\right\|^{2}_{W_{t+1}} -\left\| e_{t}\right\|^{2}_{ \Omega^{+}_{t}} -\left\|v_t\right\|^{2}}{2} \right)
\end{aligned}\end{equation*}
where $W_{t+1}:=\Omega^{+}_{t+1}+\lambda^{-1}_t  H_t^{\top} H_t$. Then, substituting (\ref{e_t+1}) we obtain \begin{equation}\label{pes_1}
\begin{aligned}
\tilde \psi_{t}\left(v_{t}| e_t\right)  \sim \exp  &\left(\frac{1}{2}\left\|( A_t- G_{t} C_t)  e_{t}+( B_t- G_{t}  D_t) v_{t}\right\|_{W_{t+1}}^{2}\right.\\
&\hspace{0.9cm}\left.-\frac{1}{2}(\left\| e_t\right\|_{ \Omega^{+}_t}-\left\|v_t\right\|^2)\right).
\end{aligned}
\end{equation}
Accordingly, the least favorable density is  \begin{equation}\label{psi_1}\tilde{\psi}_{t}\left(v_{t} | e_{t}\right) \sim \mathcal{N}(F_{t}  e_{t}, {K}_{v_{t}})\end{equation}
where \begin{align*} {K}_{v_{t}}&=\left(I_{r_t}-( B_t- G_{t}  D_t)^{\top} W_{t+1}( B_t-G_{t} D_t)\right)^{-1}\nonumber\\
 F_{t}&={K}_{v_{t}}( B_t- G_{t} D_t)^{\top} W_{t+1}(A_t- G_{t}  C_t) .\end{align*}
Then, in view of (\ref{pes_1}) and (\ref{psi_1}), it is not difficult to see that $\Omega_t^{+}$ satisfies the backward recursion:
\begin{equation*} \Omega_{t}^{+}=( A_t- G_{t} C_t)^{\top} W_{t+1}( A_t- G_{t}  C_t)+ F_{t}^{\top} {K}_{v_{t}}^{-1} F_{t}.\end{equation*}
where $\Omega^+_{T+1}=0$.

In view of (\ref{psi_1}), the least favorable noise admits the following decomposition:
$$v_t= F_t e_t+L_t \epsilon_{t}$$
where $  L_t$ is a square root matrix of $K_{v_t}$ and $\epsilon_{t}$ is normalized WGN. Therefore, taking the state space vector $$\xi_{t} \triangleq\left[\begin{array}{l}
x_{t} \\
e_{t}
\end{array}\right],$$
we obtain the corresponding least favorable model as:
\begin{equation*}\begin{aligned}
\xi_{t+1} &=\tilde{A}_{t} \xi_{t}+\tilde{B}_{t} \epsilon_{t} \\
y_{t} &=\tilde{C}_{t} \xi_{t}+\tilde{D}_{t} \epsilon_{t}
\end{aligned}\end{equation*}
where
\begin{align*}\begin{aligned}
&\tilde{A}_{t}:=\left[\begin{array}{cc}
 A_t & B_t F_{t} \\
0 &  A_t- G_{t}  C_t+\left( B_t- G_{t}  D_t\right) F_{t}
\end{array}\right]\\
&\tilde {B}_{t}:=\left[\begin{array}{c}
 B_t \\
 B_t-G_{t} D_t
\end{array}\right]L_t\\
&\tilde {C}_{t}:=\left[\begin{array}{lll}
C_t & D_t  F_{t}
\end{array}\right], \quad \tilde {D}_{t}:= D_t L_{t}.
\end{aligned}\end{align*}
It is worth noting that the least favorable model is generated by a backward recursion implementation.

Finally, it remains to evaluate the performance of the proposed robust Kalman filter and the standard Kalman filter under the least favorable model. Let $G^{\prime}_t$ be the gain of an arbitrary filter of the form
\begin{equation*}\hat{x}_{t+1}^{\prime}=A_{t} \hat{x}_{t}^{\prime}+G_{t}^{\prime}\left(y_{t}-C_{t} \hat{x}_{t}^{\prime}\right).\end{equation*}
Then, we define the corresponding prediction error $e_{t}^{\prime}=x_{t}-\hat{x}_{t}^{\prime}$. It is not difficult to see that
\begin{equation*}
\begin{aligned}
e_{t+1}^{\prime}=\left(A_t-\right.& \left.G_t^{\prime} C_t \right) e^{\prime}_t  \\
+& \left( B_t- G^{\prime}_{t} D_t \right) F_t e_t + \left( B_t- G^{\prime}_{t}  D \right) L_t \epsilon_{t}.
\end{aligned}
\end{equation*}
Taking into account (\ref{e_t+1}), we obtain
\begin{equation*}
\begin{aligned} &\left[\begin{array}{c}
e_{t+1}^{\prime} \\
e_{t+1}
\end{array}\right]=\left(\tilde{A}_{t}-\left[\begin{array}{c}
G_{t}^{\prime} \\
0
\end{array}\right]  \tilde{C}_{t}\right)\left[\begin{array}{c}
e_{t}^{\prime} \\
e_{t}
\end{array}\right]\\
&\hspace{2.1cm}+\left(\tilde{B}_{t}-\left[\begin{array}{c}
G_{t}^{\prime} \\
0
\end{array}\right] \tilde{D}_{t}\right) \epsilon_{t}.
\end{aligned}
\end{equation*}
Then, it is not difficult to see that the covariance matrix of the augmented error is given by the Lyapunov equation
\begin{equation*}\begin{aligned}
\Pi_{t+1}&=\left(\tilde{A}_{t}-\left[\begin{array}{c}
 G_{t}^{\prime} \\
0
\end{array}\right] \tilde{C}_{t}\right) \Pi_{t}\left(\tilde{A}_{t}-\left[\begin{array}{c}
 G_{t}^{\prime} \\
0
\end{array}\right] \tilde{C}_{t}\right)^{\top} \\
&\hspace{0.5cm}+\left(\tilde{B}_{t}-\left[\begin{array}{c}
 G_{t}^{\prime} \\
0
\end{array}\right] \tilde {D}_{t}\right)\left(\tilde{B}_{t}-\left[\begin{array}{c}
 G_{t}^{\prime} \\
0
\end{array}\right] \tilde {D}_{t}\right)^{\top}
\end{aligned}
\end{equation*}
where its initial value is given by $\Pi_0=\mathbbm{1}_2\otimes \tilde P_0$ and $\mathbbm{1}_2$ is the $2\times 2$ matrix whose entries are 1. Accordingly, the least favorable covariance matrix of the estimation error is the submatrix of $\Pi_t$ in position $(1,1)$.
\section{Convergence analysis} \label{sec_6}
In this section we study the convergence of the robust Kalman filter outlined in Algorithm \ref{code:recentEnd} under the assumption that the nominal model (\ref{s_space})
has constant parameters, that is $A_t=A$, $B_t=B$, $C_t=C$ and $D_t=D$, and the tolerance is constant, i.e.  $c_t=c$.
Without loss of generality we assume that $BD^{\top}=0$. Otherwise, we can rewrite the filter in Algorithm \ref{code:recentEnd} with $\tilde A= A-BD^{\top}(DD^{\top})^{-1}C$, $\tilde B$ such that  $\tilde B\tilde B^{\top}=B(I-D^{\top}(DD^{\top})^{-1}D)B^{\top}$, $\tilde C=C$ and $\tilde D=D$. In this way $\tilde B \tilde D^{\top}=0$. The robust filter converges if and only if the least favorable {\em a posteriori} covariance matrix $\tilde P_t$ of $x_t$ given $Y_{t-1}$ converges as $t\rightarrow \infty$. In view of (\ref{def_P_tilde}),
$\tilde P_t$ converges if and only if $P_t$ converges. Accordingly, the robust filter converges if and only if the iteration
\al{ \label{iteration}P_{t+1}=r_{c}(P_t),\;\; P_0\in \overline{\mathcal Q^n_+}} converges. The mapping $r_{c}$ is defined as follows
\al{ r_{c}(P_t)&=A\tilde P _tA^{\top}-A\tilde P_t C ^{\top}(C\tilde P_t C^{\top}+R)^{-1}C\tilde P_tA^{\top}+Q, \nn}
where $R:=DD^{\top}$, $Q:=B B^{\top}$, $\tilde P_{t+1}$  has been defined in (\ref{def_P_tilde})
and $\theta_t=\lambda_t^{-1}$ is such that $\gamma(P_{t+1}, \theta_t^{-1})=c$ holds.
Note that, $r_{c}(\cdot)$ is a mapping of $\overline{Q_+^n}$.
\begin{proposition} \label{prop1_conv}Let $\bar P_t$, $t\geq 0$, be the sequence generated by the Riccati equation
\al{\label{ricc_classic}\bar P_{t+1}=A\bar P _tA^{\top}-A\bar P_t C ^{\top}(C\bar P_t C^{\top}+R)^{-1}C\bar P_tA^{\top}+Q}
with $\bar P_0=P_0$. Then, the sequence $P_t$, $t\geq 0$, generated by (\ref{iteration}) is such that
\al{\label{cond_prop_ARE}\mathrm{Im}(P_t)=\mathrm{Im}(\bar P_t).} \end{proposition}
\IEEEproof We prove the first claim by induction. Clearly, condition (\ref{cond_prop_ARE}) holds for $t=0$. We assume that $\mathrm{Im}(P_t)=\mathrm{Im}(\bar P_t)$. Let $P_t=U_tD_tU_t^{\top}$ and
$\bar P_t=U_t\bar D_tU_t^{\top}$ be the corresponding reduced singular value decompositions with $D_t, \bar D_t>0$ and the columns of $U_t$ form an orthonormal basis of $\mathrm{Im}(P_t)$. Substituting such decomposition of $\bar P_t$ on the right hand side of (\ref{ricc_classic}), we obtain
\al{\bar P_{t+1}&=\tilde A_t \bar D_t \tilde A_t^{\top}-\tilde A_t \bar D_t\tilde C_t^{\top} (\tilde C_t \bar D_t \tilde C_t^{\top}+R)^{-1}\tilde C_t \bar D_t \tilde A_t^{ \top}+Q\nn\\
&=\tilde A_t (\bar D_t^{-1}+ \tilde C_t^{\top} R^{-1}\tilde C_t)^{-1}  \tilde A_t^{\top}+Q\nn}
where $\tilde A_t=AU_t$ and $\tilde C_t=C U_t$. Since $\bar D_t$ is positive definite, $v\in \mathrm{ker}(\bar P_{t+1})$ if and only if
\al{\label{cond_im}v\in \mathrm{ker}(\tilde A_t) \hbox{ and } v\in \mathrm{ker}(Q).}
In a similar way, substituting the decomposition of $P_t$ in (\ref{iteration}), we obtain
\al{P_{t+1}=\tilde A_t (D_t^{-1}-\theta_t I+\tilde C_t^{\top} R^{-1}\tilde C_t)^{-1}  \tilde A_t^{\top}+Q\nn} where we exploited the fact that $H_t^{\top}H_t=U_tU_t^{\top}$. Since $D_t^{-1}-\theta_tI$ is positive definite, $v\in \mathrm{ker}(\bar P_{t+1})$ if and only if condition (\ref{cond_im}) holds. We conclude that $\mathrm{Im}( P_{t+1})=\mathrm{Im}(\bar P_{t+1})$.
\qed\\

If the pair $(A,B)$ is stabilizable and the pair $(A,C)$ is detectable,  we have that $\bar P_t$ converges to a unique solution, say $\bar P_\infty$, for any arbitrary initial condition $\bar P_0$, see \cite{lancaster1995algebraic}. Let $U$ and $V$ be two matrices whose columns form an orthonormal basis for $\mathrm{Im}(P_\infty)$ and its orthogonal complement in $\mathbb R^n$, respectively; from now on we assume both $U$ and $V$ are fixed. By Proposition \ref{prop1_conv} it follows that
there exists $\bar t \in \mathbb N$  such that
\al{\label{red_seq}P_t^R:=U^{\top} P_t U>0, \; \; \forall \, t\geq \bar t.}
Moreover, if $P_t^R$ converges  and $V^\top P_t V\rightarrow 0$ then $P_t$ converges. Accordingly,  in what follows we prove that the sequence $P^R_t$ converges and $V^\top P_t V\rightarrow 0$ for any arbitrary initial condition.

\begin{lemma} \label{lemma_conv} Let  $P_t$, $t\geq 0$, be the sequence generated by (\ref{iteration}). We assume that $(A,B)$ and $(A,C)$ are stabilizable and detectable, respectively. If $P_0$ is such that $\mathrm{Im}(P_0)\subseteq \mathrm{Im}(U)$, then $\mathrm{Im}(P_t)\subseteq \mathrm{Im}(U)$ for any $t>0$.
Moreover,  the subsystem $(U^{\top} AU,U^{\top} B,CU,D)$ is reachable.
\end{lemma}
\IEEEproof In view of Proposition \ref{prop1_conv}, we can prove the claim using $\bar P_t$. Let $(\check A,\check B,\check C,\check D )$ be the reachability standard form of $(A, B,C,D)$ with $\check A=TA T^{-1}$, $\check B=TB$, $\check C= CT^{-1}$ and $T$ is the transformation matrix. Hence, we have
\al{\label{ABCD_basis}&\check A=\left[\begin{array}{cc}A_{11} & A_{12}  \\0 & A_{22} \end{array}\right], \; \; \check B=\left[\begin{array}{c}B_1 \\ 0\end{array}\right], \nn\\ & \check C=\left[\begin{array}{cc}C_1 & C_2\end{array}\right], \; \; \check D=D} and the pair $(A_{11},B_1)$ is reachable.  Let $\check P_t:=T\bar P_t T^{\top}$, then it is not difficult to see that $\check P_t$ obeys the following recursion
\al{\label{ricc_change}\check P_{t+1}= \check A\check{{P}}_t\check A^{\top}- \check A\check{P}_t \check C ^{\top}(\check C \check{ P}_t\check C^{\top}+\check R)^{-1}\check C\check{ P}_t\check A^{\top}+\check Q} with $\check R=R$, $\check Q=TQT^{\top}$. In plain words, (\ref{ricc_change}) is (\ref{ricc_classic}) in the new coordinates. Accordingly, we can prove the claim with respect to the new coordinates. Since the system is stabilizable and detectable,  the iteration in (\ref{ricc_change}) converges to a unique solution $\check P_\infty$ which solves the corresponding algebraic equation. By direct computation it is not difficult to see that \al{\check P_\infty=\left[\begin{array}{ccc} \check P_\infty^R& 0 \\0 & 0\end{array}\right]\nn} and $ \check P_\infty^R$ is the solution to the algebraic Riccati equation corresponding to the system $(A_{11},B_1,C_1,D)=(U^{\top}AU,U^{\top}B,CU,D)$ with $U=[\,  I_l\; 0\,]^{\top}$ and $l$ is the { dimension} of $\check P_\infty^R$. Since this subsystem is reachable, we have that $\check P_\infty^R$ is positive definite and thus
$$\mathrm{Im}(U)=\mathrm{Im}(\check P_\infty)=\mathrm{Im}\left(\left[\begin{array}{cc}I_l & 0 \\0 & 0\end{array}\right]\right).$$
Finally, by direct computation  it is not difficult to see that
\al{\check P_t=\left[\begin{array}{ccc} \star & 0 \\0 & 0\end{array}\right]\subseteq\mathrm{Im}(U)\; \implies\; \check P_{t+1}=\left[\begin{array}{ccc} \star & 0 \\0 & 0\end{array}\right]\subseteq \mathrm{Im}(U)\nn}
which concludes the proof.\qed

\begin{proposition}\label{prop_conv_subspace}
Consider the sequence generated by (\ref{iteration}) with $\mathrm{Im}(P_0)\subseteq \mathrm{Im}(U)$. Assume that $(A,B)$ is stabilizable and the reachable subsystem is also observable. Then, there always exists $c_{MAX}>0$ such that if $0<c<c_{MAX}$ then the sequence in (\ref{red_seq}) converges to a unique point $P_\infty^R$.
\end{proposition}
\IEEEproof  By Lemma \ref{lemma_conv} we have that
\al{\label{red_seq2}P_t=UP^R_tU^{\top}, \; \; t\geq 0}
and by (\ref{red_seq}) we have that $\mathrm{Im}(P_t)=\mathrm{Im}(U)$ for $t\geq \bar t$.
In view of (\ref{iteration}) and (\ref{red_seq2}), for $t\geq \bar t$ we have
\al{\label{iter_red}P^R_{t+1}=U^{\top}r_c(P_t)U=U^{\top}r_c(UP^R_tU^{\top})U=r^R_c(P^R_t)}
where
\al { r^R_c(P_t^R)&:= A_R((P_t^R)^{-1}-\theta_t I+C_R^{\top}R^{-1}C_R)^{-1} A_R^{\top}+Q_R\nn\\ A_R&:=U^{\top} A U\nn\\
Q_R&:=U^{\top} BB^{\top} U\nn\\
C_R&:= CU\nn
}  and we exploited the fact that $\mathrm{Im}(H_t^{\top})=\mathrm{Im}(P_t)=\mathrm{Im}(U)$, for $t\geq \bar t$, and thus
\al{U^{\top}\tilde P_t U &=U^{\top}(P_t^{+}-\theta_t H_t^{\top} H_t)^{+}U\nn\\&={ ((P_t^R)^{-1}-\theta_t I_l)^{-1}}.\nn} Next, we consider the subsequence $P^{Rd}_k:= P^R_{kN}$ with $k\geq \lceil \bar t/N\rceil$ and $N\geq n$. Then, from (\ref{iter_red}) we have that
\al{\label{iter_red2} P_{k+1}^{Rd}=r^{Rd}_c(P_k^{Rd})}
where $r^{Rd}_c(\cdot)$ is the $N$-fold composition of $r^{R}_c(\cdot)$ and it is defined as follows
\al{\label{iter_red2bis} r^{Rd}_c(P^{Rd})= \mathbf A_N[(P^{Rd})^{-1}+\mathbf R_N]^{-1}\mathbf A_N+\mathbf Q_N}
where $\mathbf A_N\in \mathbb R^{l \times l}$ and whose precise definition (not useful here) can be found in \cite{ZORZI_CONTRACTION_CDC,CONVTAU}; $\mathbf Q_N\in \overline{\mathcal{Q}^{l}_+}$ and $\mathbf R_N\in \overline{\mathcal{Q}^l_+}$ are defined as follows:
\al{
\mathbf Q_N &:={\cal R}_N{\cal P}_N{\cal R}_N^{\top}
\nn\\
\mathbf R_N & :=  {\cal O}_N^{\top}({\cal D}_N {\cal D}_N^{\top}+{\cal H}_N{\cal H}_N^{\top})^{-1}{\cal O}_N + {\cal J}_N^{\top} S_N^{-1} {\cal J}_N
\nn\\
S_N &:= - \Theta_N^{-1}  + {\cal L}_N(I+{\cal H}_N^{\top} ({\cal D}_N {\cal D}_N^{\top})^{-1}{\cal
H}_N)^{-1}{\cal L}_N^{\top} \nn\\
 {\cal P}_N& := [ I + {\cal H}_N^{\top} ({\cal D}_N {\cal D}_N^{\top})^{-1} {\cal H}_N -  {\cal L}_N^{\top} \Theta_N{\cal
L}_N ]^{-1} \nn\\
 {\cal J}_{N}& :=  {\cal O}_{N}^{R}-{\cal L}_{N} {\cal H}_{N}^{\top}[{\cal D}_{N} {\cal D}_{N}^{\top}+{\cal H}_{N} {\cal H}_{N}^{\top}]^{-1} {\cal O}_{N} \nn  \\
\Theta_N&:= \mathrm{diag}\left(\tilde\lambda_{N}^{-1}, \tilde \lambda_{N-1}^{-1},\ldots , \tilde\lambda_{1}^{-1}\right)\otimes I_l \nn\\
 {\cal R}_N & := \bmat {cccc}   Q_R^{1/2} &  A_R Q_R^{1/2}
& \ldots &  A_R^{N-1} Q_R^{1/2} \emat \nn\\ {\cal O}_N & := \bmat {cccc}
 (C_RA_R^{N-1})^{\top} & \ldots &  (C_R A_R)^{\top} &  C_R^{\top}\emat^{\top} \nn \\
{\cal O}_N^R & :=  \bmat {cccc}
 (A_R^{N-1})^{\top} & \ldots &  A_R^{\top} & I
\emat^{\top}  \nn\\
{\cal D}_N & := I_N \otimes  D\nn\\
{\cal H}_N & :=\mathrm{Tp} \left(
0 \; H_1,H_2 , \cdots, H_{N-2} , H_{N-1}
\right)\nn\\
{\cal L}_N & :=\mathrm{Tp} \left(
0  , L_1 , L_2  , \cdots  , L_{N-2}  ,L_{N-1}
\right)\nn\\
H_t & := \left \{ \begin{array} {cc}
 C_R A_R^{t-1} Q_R^{1/2} & t \geq 1 \\
0 & \mbox{otherwise }
\end{array} \right. \nn\\
  L_t & := \left \{ \begin{array} {cc}
A_R^{t-1} Q_R^{1/2} & t \geq 1 \\
0 & \mbox{otherwise}
\end{array} \right.\nn}
where $Q_R^{1/2}$ is a square root matrix of $Q_R$.
Then, the parameters $\tilde \lambda_l$, $l=1\ldots N$, are defined as follows. If $P^{Rd}=P^{Rd}_k$, then these parameters are given by
\al{\gamma(P_{kN+l-1}^{R}, \tilde \lambda_l)=c, \; \; l=1\ldots N.}
Since the system $(A_R,Q_R^{1/2},C_R,D)$ is reachable, by Lemma \ref{lemma_conv}, and observable, by assumption, then \cite[Proposition 5.3]{ZORZI_CONTRACTION_CDC} guarantees that there exists $c_{MAX}>0$ such that the mapping $r_c^{Rd}(\cdot)$ is a strict contraction for $0<c<c_{MAX}$ in the metric space $(\mathcal{Q}^l_+, \delta)$ where $\delta(X,Y)$ is the Thompson part metric, \cite{GAUBERT_2012,LEE_LIM_2008,LAWSON_LIM_2006}:
\al{\delta(X,Y):=\ln\max \{\sigma_{max}(Y^{-1}X),\sigma_{max}(X^{-1}Y)\}.}
Thus, there exists $0<\eta_N<1$ such that
\al{\label{contraction_rate}\delta(r_c^{Rd}(X),r_c^{Rd}(X))\leq \eta_N \delta(X,Y), \; \; \forall \, X, Y\in \mathcal Q_+^l.} It is worth noticing that $\mathbf Q_N$ and $\mathbf R_N$ in (\ref{iter_red2bis}) depend on $P^{Rd}$; however, $\eta_N$ does not. Since this metric space is complete, we can exploit the Banach fixed point theorem. Therefore, the iteration in (\ref{iter_red2}) converges to a unique point $P_\infty^R\in\mathcal Q^l_+$ for any $P_0\geq 0$ such that $\mathrm{Im}(P_0)\subseteq\mathrm{Im}(U)$. Moreover, since  $r_c^R(\cdot)$ is a mapping in $\mathcal Q_+^l$, it follows that (\ref{iter_red}) converges to $P_\infty^R$.  \qed

The next corollary is a consequence of Lemma \ref{lemma_conv} and Proposition \ref{prop_conv_subspace}.
\begin{corollary}
Consider the sequence generated by (\ref{iteration}) with $\mathrm{Im}(P_0)\subseteq \mathrm{Im}(U)$. Assume that $(A,B)$ is stabilizable and the reachable subsystem is also observable. Then, there always exists $c_{MAX}>0$ such that if $0<c<c_{MAX}$ then the sequence in (\ref{iteration}) converges to $P_\infty:=UP_\infty^RU^{\top} $.
\end{corollary}

\begin{remark} The upper bound $c_{MAX}$ for the tolerance can be computed explicitly. Following the same reasonings in \cite{ZORZI_CONTRACTION_CDC,CONVTAU} we have $c_{MAX}=\gamma(\bar P^R_q, \phi^{-1}_N)$, if $\phi^{-1}_N > \sigma_{max}(\bar P^R_q)$; otherwise $c_{MAX}=\infty$. Here: $\bar P_q^R$ is the $q$-th element of the sequence generated by the standard Riccati iteration for the reachable subsystem $(A_R,Q_R^{1/2},C_R,D)$ with initial condition $P_0^R=0$; $\phi_N\in(0,1/\sigma_{max}({\cal L}_N(I+{\cal H}_N^{\top} ({\cal D}_N{\cal D}_N^{\top} )^{-1}{\cal H}_N)^{-1}{\cal L}_N)]$ is the maximum value for which the matrix
\al{&{\cal O}_N^{\top} ({\cal D}_N {\cal D}_N^{\top} +{\cal H}_N{\cal H}_N^{\top} )^{-1}{\cal O}_N\nn\\ &\hspace{0.1cm} + {\cal J}_N^{\top}  \{- \phi_N^{-1}I+ {\cal L}_N(I+{\cal H}_N^{\top}  ({\cal D}_N {\cal D}_N^{\top} )^{-1}{\cal
H}_N)^{-1}{\cal L}_N^{\top}  \}^{ -1}{\cal J}_N\nn}
is positive definite. Clearly, the user has two degrees of freedom for computing $c_{MAX}$ that is $N$ and $q$. The larger $q$ is the larger $c_{MAX}$ is, while no specific properties in terms of $N$ have been found.\end{remark}

Proposition \ref{prop_conv_subspace} hinges on the fact that the image of $P_t$ is fixed for $t\geq \bar t$ and thus the convergence is proved in the corresponding subspace. However, in the general case, i.e. when $\mathrm{Im}(P_0)\not \subseteq\mathrm{Im}(U)$, the image of $P_t$ is not fixed because $V^\top P_t V$ could be different from the null matrix; thus, such a case needs to be addressed carefully.

\begin{lemma} \label{lemma_null_space}Consider the sequence generated by (\ref{iteration}) with $P_0\in\overline{\mathcal Q_+^n}$. Assume that $(A,B)$ is stabilizable and $(A,C)$ is detectable. Then, there always exists $c_{MAX}>0$ such that if $0<c<c_{MAX}$ then
$$ V^\top P_t V\rightarrow 0.$$
\end{lemma}
\IEEEproof We prove the claim by performing the same change of basis used in (\ref{ABCD_basis}). Let $\breve{P}_t=T P_t T^\top$, then we have
\al{\label{eqRicbreve}&\breve P_{t+1}= \check A\tilde{\breve P}_t\check A^{ \top}- \check A\tilde{\breve P}_t \check C ^{ \top}(\check C \tilde{\breve P}_t\check C^{ \top}+\check R)^{-1}\check C\tilde{\breve P}_t\check A^{ \top}+\check Q\nn\\
&\tilde{\breve P}_{t+1}=(\breve{P}_{t+1}^+-\theta_t \breve H_t^\top \breve H_t)^+} where $\breve H_t= H_t T^\top$. We partition $\breve P_t$ and $\tilde {\breve P}_t$ conformably with (\ref{ABCD_basis}):
$$\breve P_t = \left[\begin{array}{cc}\breve P_{11,t} & \breve P_{12,t}  \\\breve P_{12,t}^\top  & \breve P_{22,t} \end{array}\right],\quad \tilde{\breve P}_t = \left[\begin{array}{cc} \tilde{\breve P}_{11,t} & \tilde{\breve P}_{12,t}  \\ \tilde{\breve P}_{12,t}^\top  & \tilde{\breve P}_{22,t} \end{array}\right].$$
Taking into account (\ref{eqRicbreve}), it is not difficult to see that
\al{\label{ineq_breve}{\breve P}_{22,t+1}\leq  A_{22} \tilde{\breve P}_{22,t}A_{22}^\top.}
Notice that $\tilde {\breve P}_{t}\geq {\breve P}_{t}$ and they have the same eigenvectors corresponding to the nonnull eigenvalues, thus by (\ref{def_gamma}) we have
\al{\label{sum_c_breve}\frac{1}{2}\sum_i -\ln \left(\frac{\sigma_i(\tilde{\breve P}_{t})}{\sigma_i({\breve P}_{t})}\right)+\frac{\sigma_i(\tilde{\breve P}_{t})}{\sigma_i({\breve P}_{t})}-1= c}
where $\sigma_{i}({\breve P}_{t})$ and $\sigma_{i}(\tilde{{\breve P}}_{t})$ denote the $i$-th largest nonnull eigenvalue of ${\breve P}_{t}$ and $\tilde{{\breve P}}_{t}$, respectively.
  Since the terms on the left hand side of (\ref{sum_c_breve}) are nonnegative, it follows that
\al{\label{sum_c_breve_single} -\ln \left(\frac{\sigma_i(\tilde{\breve P}_{t})}{\sigma_i({\breve P}_{t})}\right)+\frac{\sigma_i(\tilde{\breve P}_{t})}{\sigma_i({\breve P}_{t})}-1\leq 2c.} By (\ref{sum_c_breve_single}) it follows that
\al{\label{cond_P_breve}\tilde{\breve P}_{t}\leq \rho {\breve P}_{t}}
where $\rho>1$ is the unique solution to the equation
\al{\label{rho_breve}-\ln \rho+\rho-1=2c. }
Substituting (\ref{cond_P_breve}) in (\ref{ineq_breve}), we obtain
\al{ {\breve P}_{22,t}\leq  \rho A_{22} {\breve P}_{22,t}A_{22}^\top.\nn} Accordingly, if $\sqrt{\rho}|\sigma_1(A_{22})|<1$ then ${\breve P}_{22,t}\rightarrow 0$ which also implies ${\breve P}_{12,t}\rightarrow 0$ because ${\breve P}_{t}\geq 0$. Since in these coordinates $V=[\,0 \; I_{n-l}\,]^\top$, it follows that $V^\top {\breve P}_{t}V\rightarrow 0$.

The existence of $c_{MAX}>0$ such that  $\sqrt{\rho}|\sigma_1(A_{22})|<1$ follows form the fact that the left hand side of (\ref{rho_breve}) is monotone increasing for $\rho>1$ and it approaches zero as $\rho\rightarrow 1$.\qed

\begin{theorem}
Consider the sequence generated by (\ref{iteration}) with $P_0\in \overline{\mathcal Q_+^n}$. Assume that $(A,B)$ is stabilizable and the reachable subsystem is also observable. Then, there always exists $c_{MAX}>0$ such that if $0<c<c_{MAX}$ then the sequence in (\ref{red_seq}) converges  to $P_\infty^R$.
\end{theorem}
\IEEEproof Consider the fixed point $P_\infty^R$ of Proposition \ref{prop_conv_subspace} and the subsequences $P^{d}_k:=P_{kN}$ and $P^{Rd}_k:=P^R_{kN}$ with $k\geq \lceil \bar t/N\rceil$, $N\geq n$ and $\bar t$ is such that (\ref{red_seq}) holds. Accordingly, $P^R_\infty$ and $P^{Rd}_k$, $k\geq \lceil \bar t /N \rceil$, are two elements of the metric space $(\mathcal Q_+^l,\delta)$.  Moreover, we have
 \al{P_{k+1}^d=r_c^d(P_k^d), \; \; P_{k+1}^{Rd}=r_c^{Rd}(P_k^d)\nn}
 where $r_c^d(\cdot),r_c^{Rd}(\cdot)$ are the $N$-fold composition of $r_c(\cdot),r_c^R(\cdot)$, respectively. Since $r_c(\cdot)$ is a continuous mapping of $\overline{\mathcal Q^n_+}$, then also $r_c^d(\cdot)$ is a continuous mapping of $\overline{\mathcal Q^n_+}$. By Lemma \ref{lemma_null_space} we have that $V^\top P_tV\rightarrow 0$ which implies $V^\top P_k^dV \rightarrow 0$ and hence
 \al{
r_c^d(P_k^d)=r_c^d(UU^{\top}  P_k^dUU^{\top} )+M_k\nn}
 with $M_k\rightarrow 0$. Accordingly,
\al{P_{k+1}^{Rd}&=U^{\top} r_c^d(P_k^{d}) U=U^{\top} r_c^{d}(UU^{\top}P_k^{d}UU^{\top}) U+N_k\nn\\ &=r_c^{Rd}(P_k^{Rd})+N_k\nn}
where $N_k=U^{\top} M_kU\in \mathcal Q^l$ and $\|N_k\|$ tends to zero as $k\rightarrow \infty$. Since the Thompson part metric is continuous with respect to the arguments, and in view of (\ref{contraction_rate}), we have that
\al{\delta(P_{k+1}^{Rd}, P_\infty^R)&=\delta(r_c^{Rd}(P_k^{Rd})+ N_k, P_\infty^R)\nn\\ &=\delta(r_c^{Rd}(P_k^{Rd}), P_\infty^R)+\varepsilon_k\nn\\ &=\delta(r_c^{Rd}(P_k^{Rd}), r_c^{Rd}(P_\infty^R))+\varepsilon_k\nn\\
&\leq \eta_N \delta(P_k^{Rd}, P_\infty^R)+\varepsilon_k\nn}
where $\varepsilon_k\rightarrow 0$ as $k\rightarrow \infty$. Accordingly,
\al{\delta(P_{k+1}^{Rd}, P_\infty^R)&\leq \eta_N \delta(P_k^{Rd}, P_\infty^R)+\varepsilon_k\nn\\
&\leq \eta_N^2 \delta(P_{k-1}^{Rd}, P_\infty^R)+\eta_N\varepsilon_{k-1}+\varepsilon_k\leq \ldots \nn\\ & \leq  \eta_N^{k+1}\delta(P_{0}^{Rd}, P_\infty^R)+\sum_{l=0}^k \eta_N^{k-l}\varepsilon_l\nn\\
 & = \eta_N^{k+1}\delta(P_{0}^{Rd}, P_\infty^R)+ h \ast\varepsilon(k) \nn}
where $h \ast\varepsilon (k)$ denotes the convolution in $k$ between the sequences $h_k=\eta_N^k$ and $\varepsilon_k$ with $k\geq 0$. Accordingly, $h \ast\varepsilon (k)$ can be understood as the output at time $k$ of a stable first order filter fed by a signal which tends to zero in the steady state. As a consequence $h \ast\varepsilon (k)\rightarrow 0$ as $k\rightarrow \infty$. Accordingly, we conclude that
\al{\delta(P_{k+1}^{Rd}, P_\infty^R)\rightarrow 0\nn}
as $k\rightarrow \infty$ and thus $P_{k}^{Rd}\rightarrow P^R_\infty$. It is not difficult to see that the same conclusion holds in the case that $P^{Rd}_k:=P_{kN+l}^R$ with $l=1\ldots N-1$. Accordingly, we have that $P^R_t\rightarrow P_\infty^R $ as $t\rightarrow \infty$. \qed

Finally, in view of the fact that $V^\top P_t V\rightarrow 0$ we have the following corollary.

\begin{corollary} Assume that $(A,B)$ is stabilizable and the reachable subsystem is also observable. Then, there always exists $c_{MAX}>0$ such that if $0<c<c_{MAX}$ then the sequence in (\ref{iteration}) converges to $P_\infty=U P_\infty^RU^{\top}$ for any initial condition $P_0\in\overline{\mathcal Q_+^n}$.
\end{corollary}

\section{Simulation results}\label{sec_7}
To evaluate the effectiveness of the robust Kalman filter proposed in Algorithm \ref{code:recentEnd}, we consider the nominal model:
\begin{equation*}
\begin{aligned}
x_{t+1} &=A x_{t}+B v_{t} \\
y_{t} &=C x_{t}+D v_{t}
\end{aligned}
\end{equation*}
where
\begin{align}
A&=\left[\begin{array}{ccc}
2 & 0.1 & 0.1 \\
0 & 0.8407 & -0.3482\\
0 & 0.3482 & 0.8407
\end{array}\right],\; \;  B=\left[\begin{array}{cc}
1 & 0\\
0 & 0\\
0 & 0
\end{array}\right], \nonumber\\
C&=\left[\begin{array}{lll}
1 & 0 & 0
\end{array}\right], \; \;  D=\left[\begin{array}{ll}
 0 & 1
\end{array}\right], \nonumber
\end{align} and $x_{0} \sim \mathcal{N} (0, \tilde P_{0})$ with
$$ \tilde P_0=\left[\begin{array}{lll}
1 & 0 & 0 \\
0 & 0 & 0 \\
0 & 0 & 1
\end{array}\right].$$
The pairs $(A,B)$ and $(A,C)$ are stabilizable and observable, respectively. However, it is worth noticing that the model is not reachable.
\begin{figure}[htb]
\centering
\includegraphics[width=0.5\textwidth]{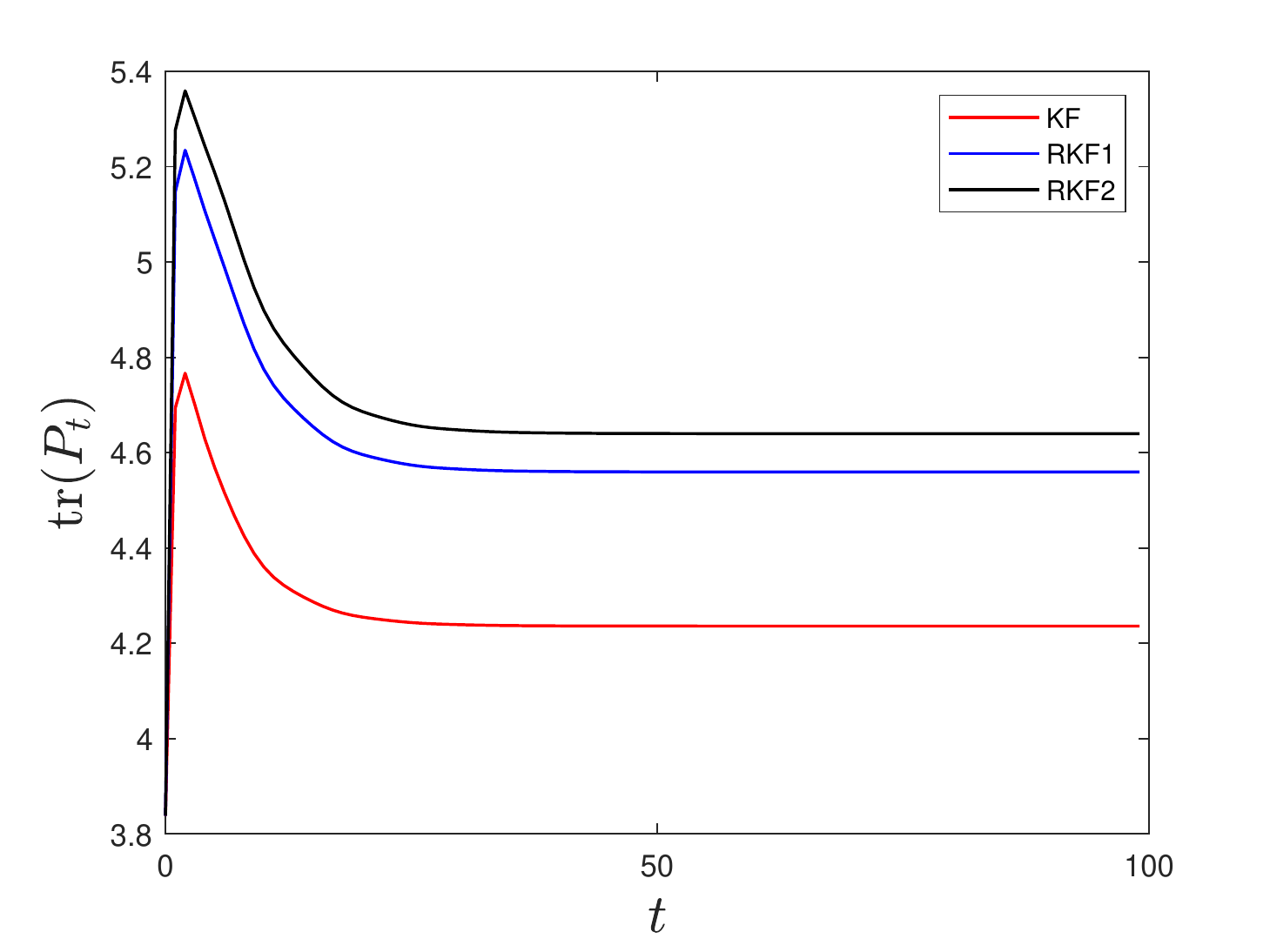}
\caption{Trace of $P_t$ for KF, RKF1 with $c= 10^{-1}$, and RKF2 with $c=2 \cdot 10^{-1}$.} \label{fig1}
\end{figure}
\begin{figure}[htb]
\centering
\includegraphics[width=0.5\textwidth]{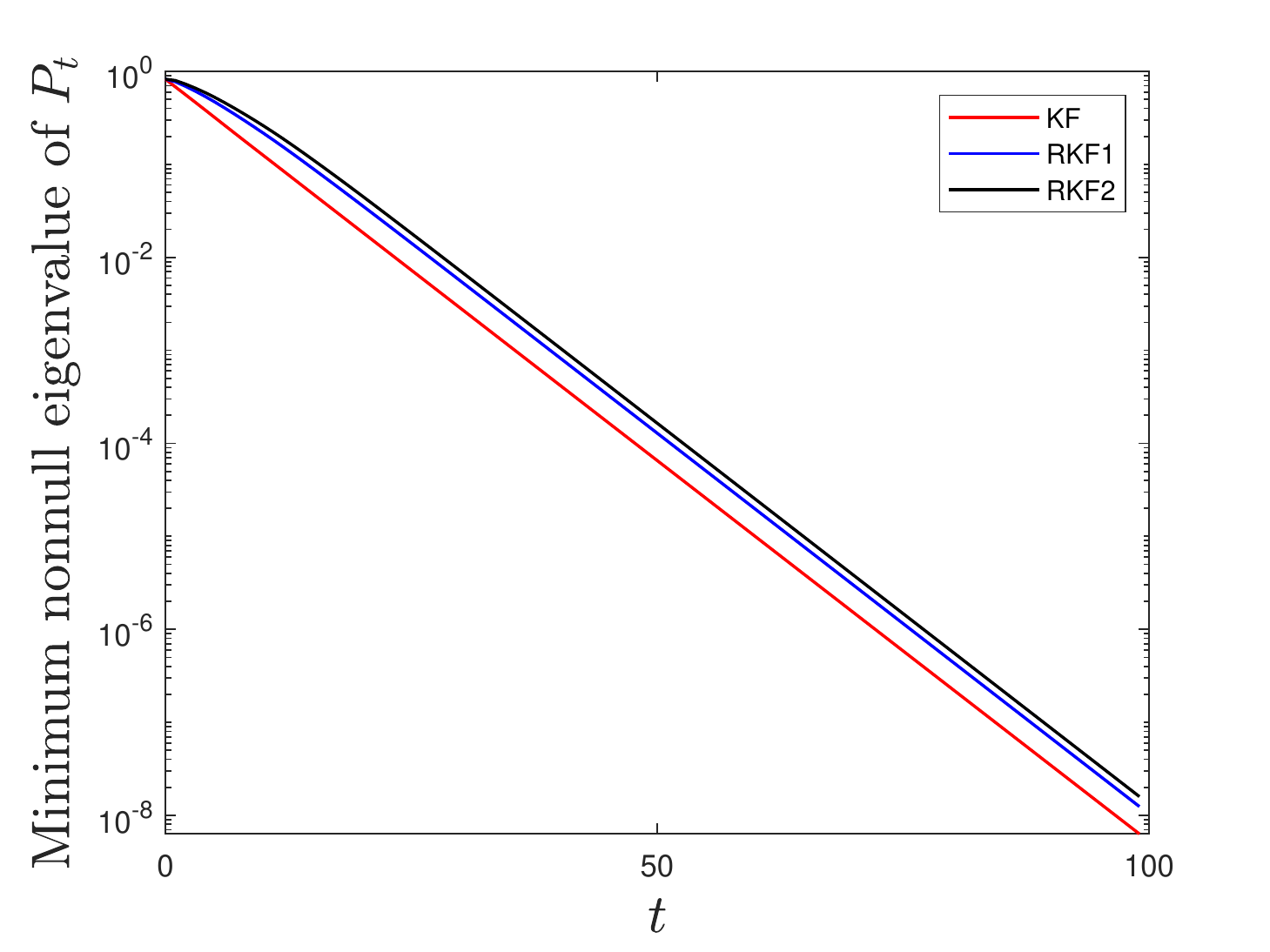}
\caption{Minimum nonnull eigenvalue of $P_t$ as a function of $t$ for KF, RKF1 with $c=10^{-1}$, and RKF2 with $c=2 \cdot 10^{-1}$.} \label{fig2}
\end{figure}
\begin{figure}[htb]
\centering
\includegraphics[width=0.5\textwidth]{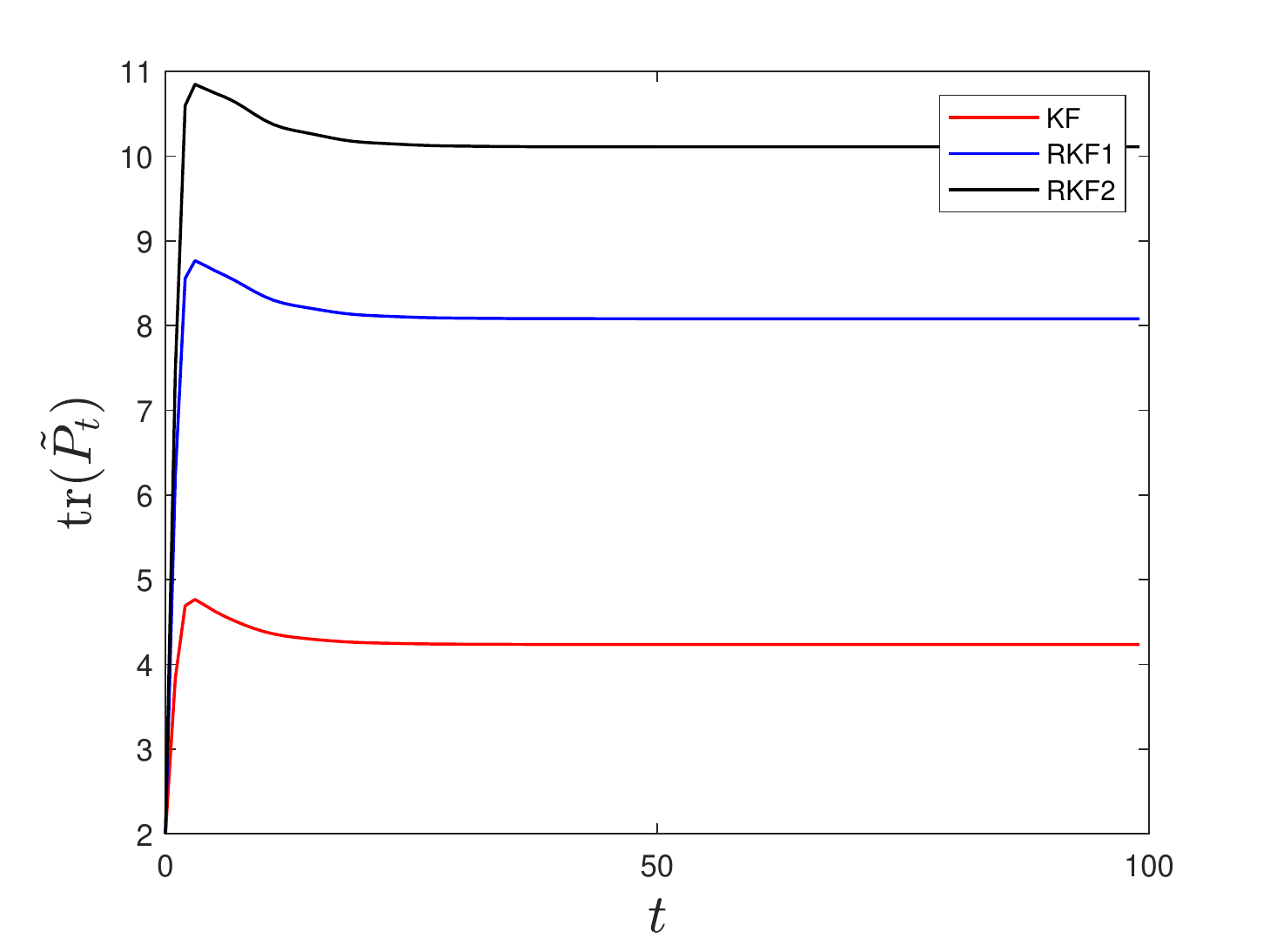}
\caption{Trace of $\tilde P_t$ for KF, RKF1 with $c= 10^{-1}$, and RKF2 with $c=2 \cdot 10^{-1}$. Recall that $\tilde P_t=P_t$ for KF.} \label{fig3}
\end{figure}

 \begin{figure}[htb]
\centering
\includegraphics[width=0.5\textwidth]{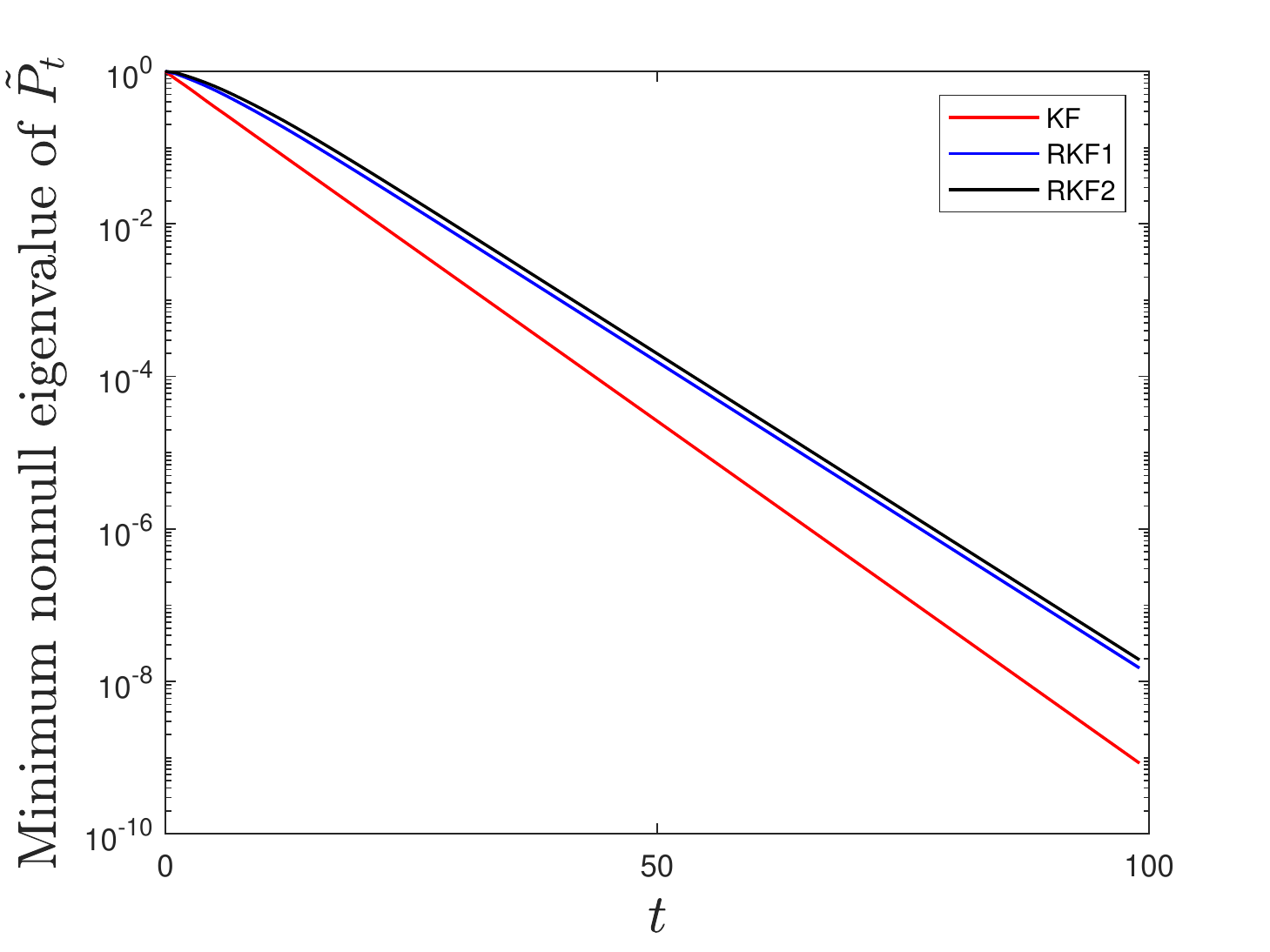}
\caption{Minimum nonnull eigenvalue of $\tilde P_t$ as a function of $t$ for KF, RKF1 with $c=10^{-1}$, and RKF2 with $c=2 \cdot 10^{-1}$.} \label{fig4}
\end{figure}

\begin{figure}[htb]
\centering
\includegraphics[width=0.5\textwidth]{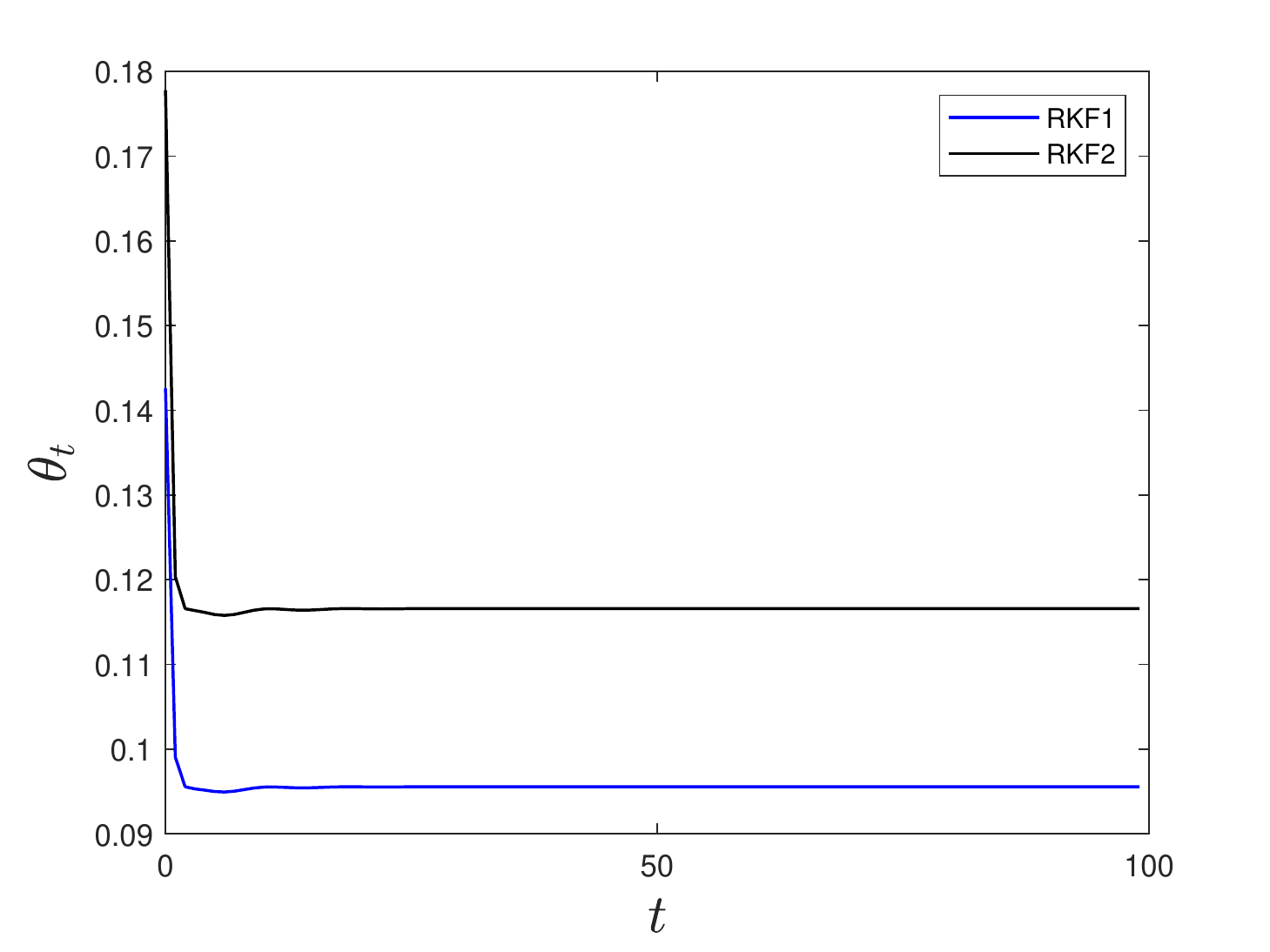}
\caption{Risk sensitivity parameter $\theta_t$ as a function of $t$ for RKF1 with $c=10^{-1}$, and RKF2 with $c=2 \cdot 10^{-1}$.} \label{fig5}
\end{figure}

In what follows, we consider the following three filters over a time interval of length $T=100$: the standard Kalman filter (KF), the robust Kalman filter (RKF) proposed in Algorithm \ref{code:recentEnd} with tolerance $c_1= 10^{-1}$ (RKF1) and with tolerance $c_2=2 \cdot 10^{-1}$ (RKF2). In Fig. \ref{fig1} we show the trace of $P_t$ generated by these three filters: it converges to a constant value in all the three cases. It is worth noting that $\mathrm{rank}(P_t)=2$ for all the algorithms. Fig. \ref{fig2} shows the minimum nonnull eigenvalue of $P_t$ as a function of the time: as time increases, its eigenvalue decreases exponentially which means that $P_t$ has numerical rank equal to 1 as time takes large values. Hence, the corresponding algebraic Riccati equation converges to a unique solution with rank equal to 1, which is in line with our expectations. It is interesting to note that the convergence rate of KF is higher than that of RKF. More precisely, the smaller $c$ is, the higher the convergence rate is. The trace of $\tilde P_t$ for KF, RKF1 and RKF2 is shown in Fig. \ref{fig3}. Clearly $\tilde P_t=P_t$ for KF. Also in this case $\mathrm{rank}(\tilde P_t)=2$ and its numerical rank is equal to 1 for  large values of $t$ in all the algorithms, see the minimum nonnull eigenvalue of $\tilde P_t$ depicted in Fig. \ref{fig4}. Finally, as shown in Fig. \ref{fig5}, the risk sensitivity parameter $\theta_t$ converges to a constant value for both RKFs. Moreover, the smaller $c$ is, the smaller $\theta_t$ is.

\begin{figure}[htb]
\centering
\includegraphics[width=0.5\textwidth]{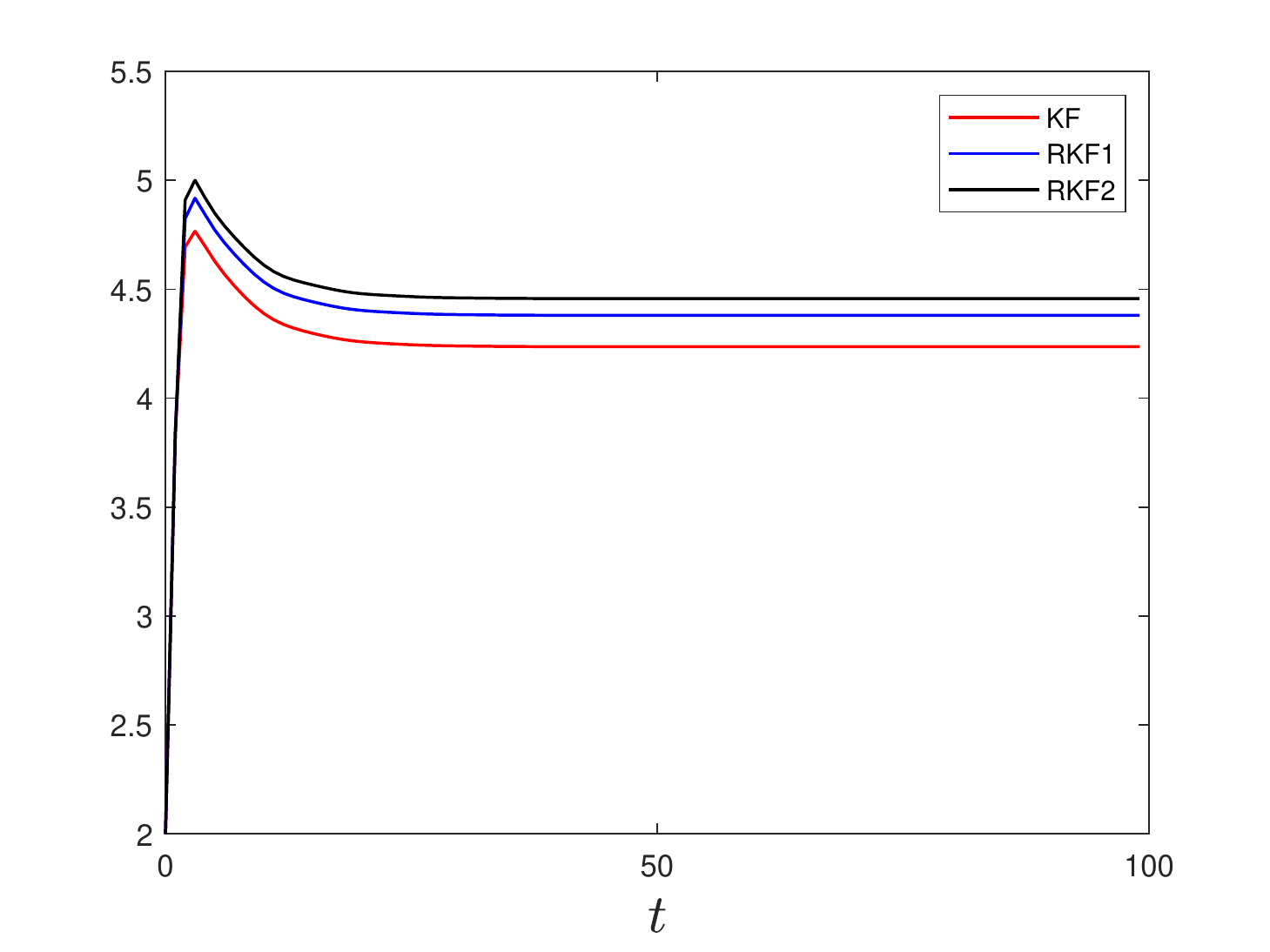}\caption{Scalar variance of the prediction error for KF, RKF1 with $c= 10^{-1}$ and RKF2 with $c= 2 \cdot 10^{-1}$ under the assumption that the actual model coincides with the nominal model.} \label{fig6}
\end{figure}

\begin{figure}[htb]
\centering
 \includegraphics[width=0.5\textwidth]{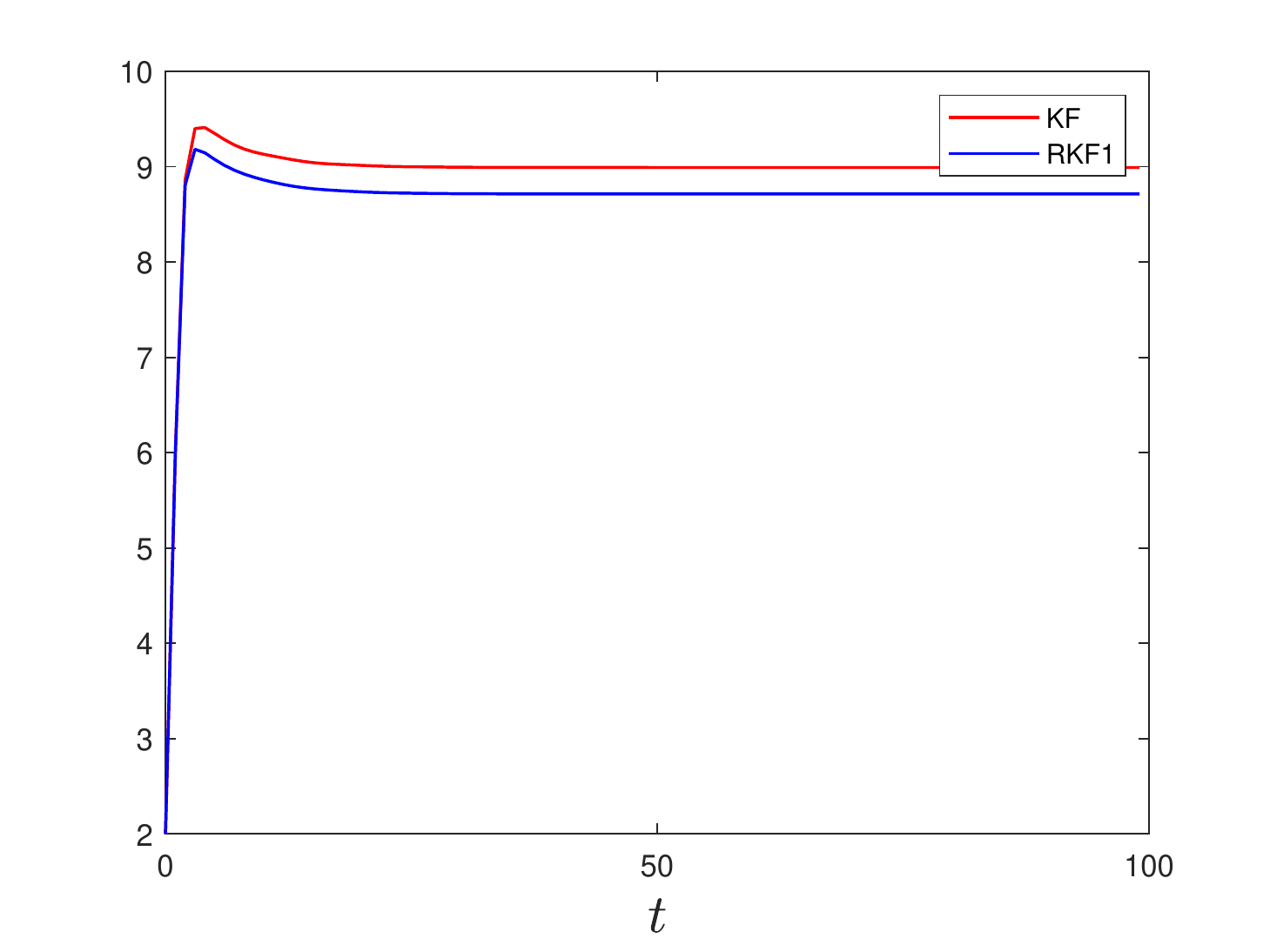}
\caption{Scalar variance of the prediction error  for KF and RKF1 with $c= 10^{-1}$ under the assumption that the actual model is the least favorable model solution to the corresponding minimax problem.} \label{fig7}
\end{figure}

\begin{figure}[htb]
\centering
 \includegraphics[width=0.5\textwidth]{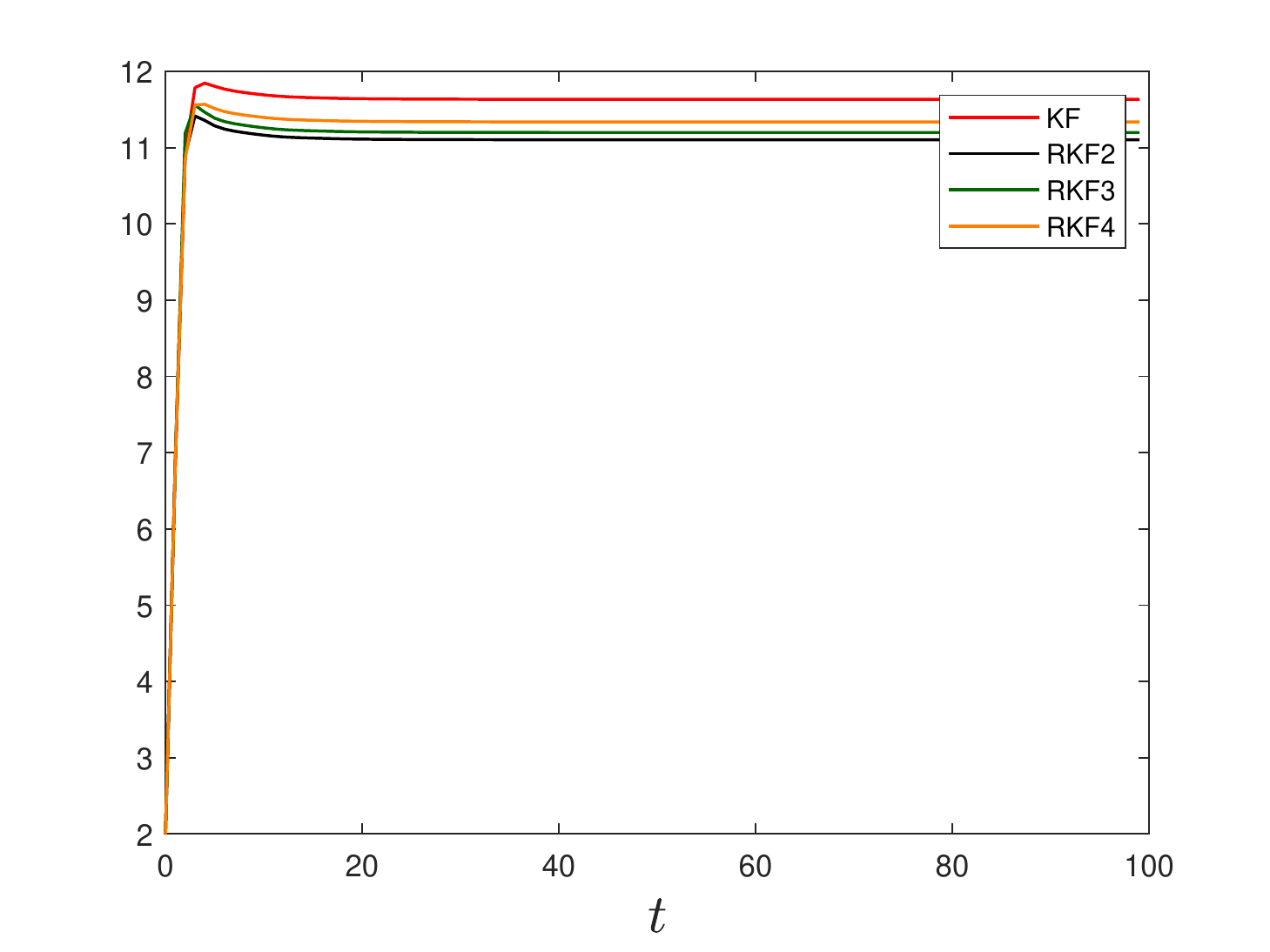}
\caption{ Scalar variance of the prediction error for KF, RKF2 with $c=2 \cdot 10^{-1}$, RKF3 with $c=1$ and RKF4 with $c= 10^{-2}$ under the assumption that the actual model is the least favorable model solution to the minimax problem with $c=2 \cdot 10^{-1}$.} \label{fig8}
\end{figure}

Fig. \ref{fig6} shows the variance of the prediction error in the case that the actual model corresponds to the nominal model. In this case, KF performs better than RKF's. On the other hand, when the actual model is the least favorable model solution to the corresponding minimax problem, KF is worse than RKF's. Fig. \ref{fig7} shows the variance of the prediction error for KF and RKF1 in the case that the actual model corresponds to the least favorable model solution to (\ref{minimax}) with $c=10^{-1}$. A similar scenario occurs in the case that we compare KF, RKF2
and the actual model is the least favorable one solution to (\ref{minimax}) with $c=2\cdot 10^{-1}$, see Fig. \ref{fig8}. The latter figure shows also the performance of the robust Kalman filters with tolerance $c=1$ (RKF3) and
$c=10^{-2}$ (RKF4). As expected, RKF3 and RKF4 are worse than RKF2, however, the former outperform KF. This means that, even though the ambiguity set is chosen  slightly large (case with RKF3) or small (case with RKF4) as compared to the actual model, RKF still outperforms KF.
Finally, it is worth noting that the least favorable model depends on tolerance $c$, and thus least favorable performance of KF depends $c$  as well, see Fig. \ref{fig7} and Fig. \ref{fig8}.

\section{Conclusions}\label{sec_8}
In this paper, we have considered a robust static estimation problem in the case that the nominal density and the actual density could be possibly degenerate. This approach is characterized by a minimax game: one player selects the least favorable density in a prescribed ambiguity set, while the estimator is designed according to this least favorable density. Then, we have extended this framework to the dynamic case where we have a state space model whose transition probability density is possibly degenerate. The solution is a Kalman-like filter whose gain is updated according to a Riccati-like iteration which evolves on the cone of the positive semidefinite matrices. We have proved that the resulting robust filter converges in the case that the nominal state space model has constant parameters and the ambiguity set has a radius $c$ sufficiently small. Moreover, we have characterized the least favorable model over a finite horizon. Finally, some numerical experiments have been presented in order to show the effectiveness of the proposed robust filter.
\section{Acknowledgments}
Shenglun Yi is partially supported by NSFC under Grant Nos. 61973036 and 61433003.

\end{document}